\documentclass[sanskrit]{article}
in \oddsidemargin 0.25true in \evensidemargin 0.25true in
\usepackage{epsfig}

\usepackage{enumerate}

\usepackage{graphicx}
\usepackage{subcaption}

\usepackage{CJKutf8}
\usepackage{devanagari}

\usepackage{amsmath}

\usepackage{lipsum}

\begin{document}

\newcommand{\arcsech}{\mbox{arcsech}}
\newcommand{\bc}{\begin{center}}
\newcommand{\ec}{\end{center}}
\newcommand{\be}{\begin{equation}}
\newcommand{\ee}{\end{equation}}
\newcommand{\bea}{\begin{eqnarray}}
\newcommand{\eea}{\end{eqnarray}}
\newcommand{\bean}{\begin{eqnarray*}}
\newcommand{\eean}{\end{eqnarray*}}
\newcommand{\bd}{\begin{description}}
\newcommand{\ed}{\end{description}}

\newcommand{\caln}{\cal N}
\newcommand{\caloi}{{\cal OI}}
\newcommand{\hatv}{{\bf v}}
\newcommand{\hatrho}{\hat{{\bf \rho}}}
\newcommand{\hatlam}{\lambda}
\newcommand{\tillam}{\tilde{\lambda}}
\newcommand{\tT}{\tilde{T}}
\newcommand{\hT}{\hat{T}}
\newcommand{\bT}{{\bf T}}
\newcommand{\bA}{{\bf A}}
\newcommand{\ep}{\epsilon}
\newcommand{\cd}{{\cal D}}
\newcommand{\hp}{\hat{T}}
\newcommand{\epf}{\epsilon_f}
\newcommand{\Nn}{N_{net}}
\newcommand{\Ns}{N_{sph}}
\newcommand{\Nd}{N_{dat}}
\newcommand{\lng}{\langle}
\newcommand{\rng}{\rangle}
\newcommand{\vep}{\epsilon}
\newcommand{\vph}{\varphi}
\newcommand{\vp}{\varphi}
\newcommand{\bfk}{{\bf k}}
\newcommand{\bfr}{{\bf r}}
\newcommand{\bfN}{{\bf N}}
\newcommand{\bfa}{{\bf a}}
\newcommand{\bfb}{{\bf b}}
\newcommand{\bfc}{{\bf c}}
\newcommand{\bfH}{{\bf H}}
\newcommand{\bfK}{{\bf K}}
\newcommand{\bfL}{{\bf L}}
\newcommand{\bfM}{{\bf M}}
\newcommand{\bfx}{{\bf x}}
\newcommand{\bC}{{\bf C}}
\newcommand{\bp}{{\bf{p}}}
\newcommand{\bfe}{{\bf e}}
\newcommand{\hbt}{\hat{\overline{T}}}
\newcommand{\uvr}{{\bf \hat{r}}}
\newcommand{\uvn}{{\bf \hat{n}}}
\newcommand{\bw}{\overline{w}}
\newcommand{\bx}{\mbox{\boldmath $\xi$}}
\newcommand{\bnu}{\mbox{\boldmath $\nu$}}
\newcommand{\my}{M_{yr}}

\newcommand{\bvt}{\begin{verbatim}}
\newcommand{\evt}{\end{verbatim}}
\newcommand{\dbar}{\vrule width 0.25truein height 1.0pt depth -0.6pt}

\renewcommand{\baselinestretch}{2.0}

{\raggedleft
{\bf\large
OFFER (One-Formula-Facilitates-Every-Relation) for Calculus }
}
\vskip 0.5cm

\renewcommand{\baselinestretch}{1.2}

{\raggedleft 
{\large SAMUEL S.P. SHEN }  \\
San Diego State University, San Diego, CA 92182, USA\\
{\tt sshen@sdsu.edu}\\
ORCID: 0000-0002-4535-4684

\vskip 0.0cm

}
\vskip 0.5cm

\thispagestyle{empty}
\renewcommand{\baselinestretch}{1.2}

\noindent   {\bf Abstract:} The large number of calculus formulas makes learning calculus heavily memory-dependent and can scare away many students. Based on Issac Newton's idea of prediction, this paper shows that calculus students need to master only one critical formula: his proposition that, given boundary data, one can predict inside and vice versa.  We can treat other formulas as consequences of the formula for this proposition.  In 1-, 2-, and 3-dimensional spaces, this formula corresponds to the fundamental theorem of calculus, Green's theorem, and the divergence theorem.  This approach is summarized as OFFER (one-formula-facilitates-every-relationship) for calculus education. This paper demonstrates the method and advantages of the OFFER pedagogy. The OFFER approach can accelerate the pace of calculus teaching, provide deeper insight into calculus,  inspire questions from students, and naturally facilitate the inclusion of real-world data, artificial intelligence (AI), and computer coding in the calculus series. Consequently, OFFER may serve as a way to modernize the calculus curricula in the era of AI and the digital economy.


\vspace{10px}

\section{What is the Main Purpose of Calculus?}
What was Issac Newton's answer to this question? What was Gottfried Wilhelm von Leibniz's answer? What are today's students' answers? What are today's instructors' answers? What are the main features of calculus applications in modern times? 

Questions like this can go on and on about the importance of calculus and why students need to learn it. The answers can be many times more than the questions. For example, calculus is a study of change, as you may have read in many calculus textbooks or hear from many instructors. However, students are puzzled by this answer because most calculus materials in the traditional engineering calculus series (Cal I, II, III) focus on integration and derivative techniques, except for the time-derivative section(s) that explicitly explain the rate of change. 

Some students may cynically think that the main purpose of calculus is to destroy them and block them from entering the field of science, technology, engineering, and mathematics (STEM). Their rationale is reasonable because most material in the conventional Cal I, II, and III series is never used, even by professional researchers, in their lifetime, such as the logarithmic derivative and half-angle substitution for finding antiderivatives. This rationale is harmful to calculus education because it demoralizes students and sometimes makes them angry. 

This cynical idea is enhanced by some instructors who treat calculus as a gatekeeper that blocks ``unqualified" or ``unworthy" students from entering the STEM fields. They might recall the ``excellent mathematical skills" of STEM students in the good old days, when only 20\% or less of high school graduates could go to college. 

Still, some instructors may think that calculus is so beautiful that young people should learn how to appreciate its beauty and power. 
If all the calculus students had this ability, the world would be very beautiful! Of course, this utopian world does not exist. 

Ironically, the instructors who emphasize beauty, power, and gatekeeping rarely have experience solving real-world problems with calculus or demonstrating its power. Even today, many of them cannot do AI-assisted computer programming for calculus research or applications using real-world data, such as measuring the surface area of a lake using Green's theorem and a GPS, let alone explaining Green's theorem and its relationship to Faraday's law in physics for electricity generation with a moving magnet. 

\begin{figure} [h]
\centering
\includegraphics[width=4.1in]{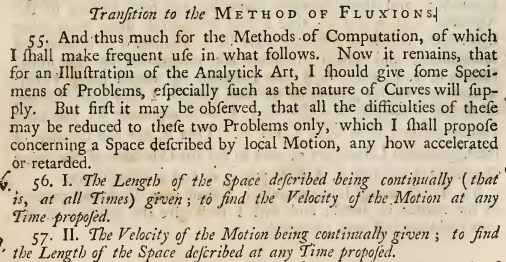}
\caption{A screenshot of Isaac Newton's statement on the purpose of calculus in his book ``{\it The Method of Fluxions and Infinite Series}" (page 19) [1].}
\label{fig01}
\end{figure}

 What did Newton say about the purpose of calculus? In his book ``The Method of Fluxions and Infinite Series; with Its Application to the Geometry and Curve-Lines," Isaac Newton considered the ``Method" as a prediction tool. 
Figure \ref{fig01} shows a screenshot from the chapter of 
``Transition to the Method of Fluxions" in his book on the method of fluxions, commented in 1736
by John Colson [1]. Newton envisioned that "all the difficulties of these may be reduced to these two problems only: Problems 56 and 57 shown in Fig. \ref{fig01} (Page 19 of [1]). 
\begin{description}
\item{Problem 56.} ~ I. The length of the Space described being continually (that is, at all Times) given, to find the Velocity of the Motion at any Time proposed.
\item{Problem 57.} ~ II. The Velocity of the Motion being continually given; to find the Length of the Space described at any Time proposed.
\end{description}
 
You may use a car's odometer and speedometer data to understand his two problems. Given the odometer record at all times, find the car's speed at any specified time.  Given the speedometer record at all times, find the distance the car has traveled from the beginning to any specified time. We can see that Newton formulated a precise statement of two prediction problems. 

In his book [1], he solved Problem 56 using the binomial series method, and Problem 57 using the antiderivative method and boundary data: 
The traveled distance equals the endpoint odometer reading minus the startpoint reading. This is a layman statement of the fundamental theorem of calculus (FTC). 

In today's calculus language, his two problems may be translated as follows:
\begin{enumerate}
\item[(A)] Given the data of the function $f(x)$ at all times $x$ continuously, find the derivative $f'(x)$ at any specified time.
\item[(B)] Given the derivative data $f'(x)$ at all times $x$ continuously, find the function data $f(x)$ at at any specified time. 
\end{enumerate}

The solution to Newton's Problem 56 is to find the derivative $f'(x)$ (or $\dot{y}$ called fluxion by Newton), 
and that to his Problem 57 is to find an integral:
\be
 \int_a^xf'(t) ~\mbox{d}t = f(x) - f(a) 
 \label{eq01ftc}
\ee
where $a$ is the starttime and $x$ is the specified time. This is an analytic expression of FTC using modern calculus notation. 

A special case is $x=b$ in the interval $[a, b]$, then Eq. (\ref{eq01ftc}) becomes
\be
 \int_a^bf'(t) ~\mbox{d}t = f(b) - f(a) = f(x)\bigg|_a^b.
 \label{eq02ftc}
\ee
This expression shows that the boundary data of the original function can be used to predict a type of derivative accumulation inside $[a, b]$. 

We argue in this article that calculus as a prediction tool is a more appropriate way of thinking for today's calculus education than the ideas presented in most calculus books on the market. We show how the single FTC formula can be used to derive all the important relationships in calculus in 1- to 3-dimensional (3D) settings: 1D corresponds to the FTC, 2D to Green's theorem, and 3D to the divergence theorem.  The mean value theorem and Taylor's theorem are consequences of the FTC prediction formula. With these, we provide an approach to one-formula calculus named OFFER (one-formula-facilitates-every-relationship). We hope that the OFFER approach can make calculus learning easier, more insightful, and deeper. Calculus learning becomes easier because (i) students do not need to memorize many formulas as in the conventional textbooks, (ii) they learn based on their common sense and hand-drawn figures rather than the traditional abstract ideas and complex notations, and (iii) they can focus only on a single fundamental idea and extend it to numerous other cases as Newton did in his book [1]. 

More importantly, in the era of the digital economy, OFFER makes it natural to include real-world data, artificial intelligence (AI), and computer coding in the calculus courses. This can empower our students to use calculus and AI to solve practical problems and genuinely appreciate the power and beauty of calculus. 

\section{Define the Integral as the Exact Sum of Staircase Risers without Approximation}



\begin{figure} [h]
\centering
\includegraphics[width=3.0 in]{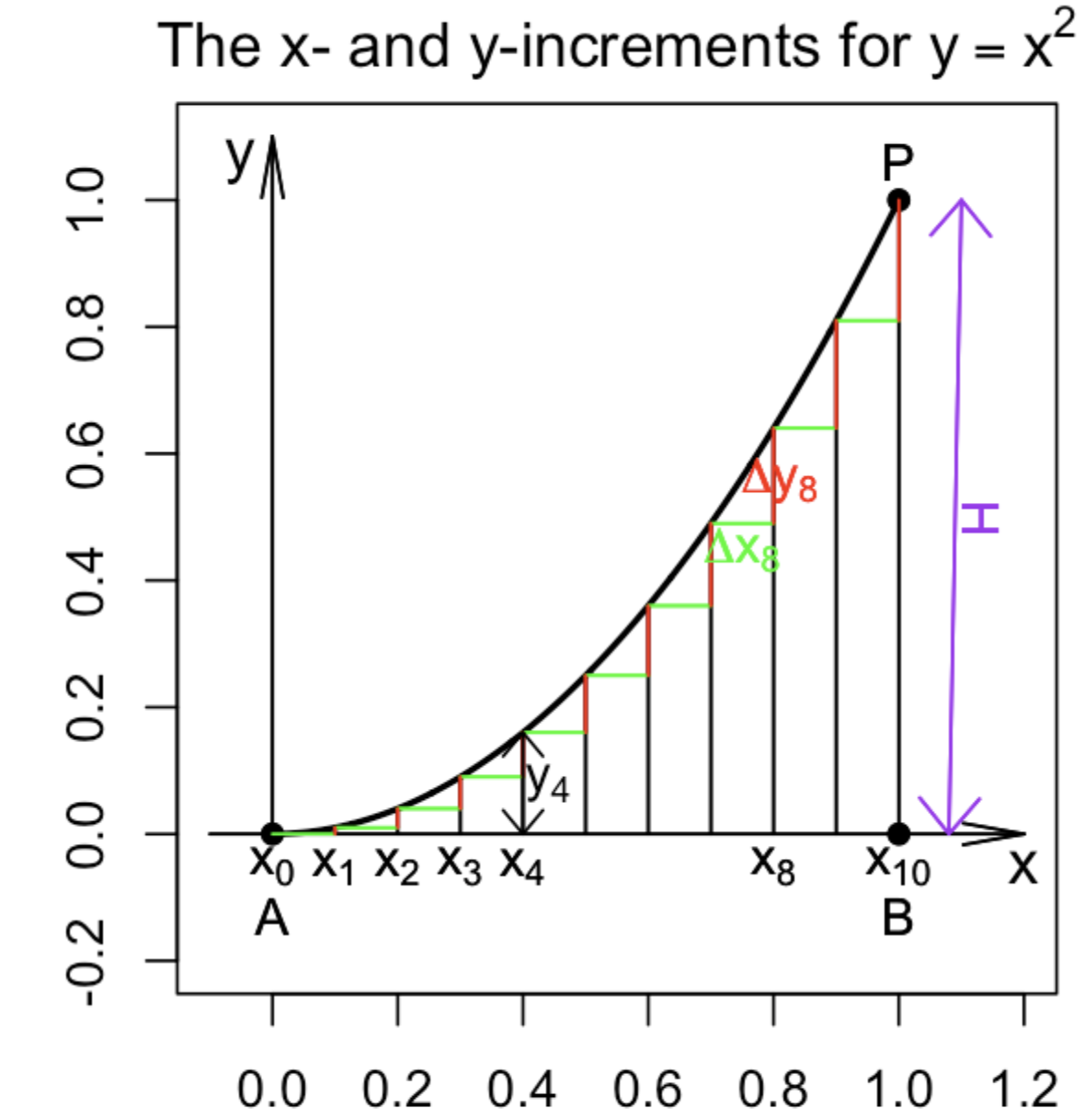}
\caption{ The staircases from point $A$ to $P$ with the total rise $H$, where $H$ is equal to the sum of all the step risers, expressed explicitly by Eq. (\ref{eq030}).
The black curve is determined by the function $y = f(x) = x^2$. }
\label{fig02}
\end{figure}

The staircase and the curve $y=f(x)$ in Fig. \ref{fig02} implies the following
\bea
H &=& f(b) - f(a)  \nonumber\\
&=&f(x_n) - f(x_0) \nonumber\\
&=& f(x_n) - f(x_{n-1}) + f(x_{n-1}) - f(x_{n-2}) + \nonumber\\
&&\cdots  + f(x_2)- f(x_1) + f(x_1) - f(x_0) \nonumber\\
&=&  \big(f(x_n) - f(x_{n-1})\big) + \big(f(x_{n-1}) - f(x_{n-2})\big) +  \nonumber\\
&& \cdots  + \big(f(x_2) - f(x_1)\big)+ \big(f(x_1) - f(x_0)\big) \nonumber\\
&=& \Delta y_n + \Delta y_{n-1} + \cdots +  \Delta y_2 +  \Delta y_1.
\label{eq030}
\eea
By introducing $ - f(x_{n-1}) + f(x_{n-1})$ etc., we explicitly express that the total height $H$ is equal to the 
sum of the rises of all the staircases. Most ten-year-olds should possess this level of common sense or intuition.

Now, we introduce a so called integral notation $\int_{D_y} \mbox{d}y$ for this summation:
\be
 \int_{D_y} \mbox{d}y = \Delta y_n + \Delta y_{n-1} + \cdots +  \Delta y_2 +  \Delta y_1,
 \label{eq0304}
\ee
where $\int$, like a lowercase $s$ in the old English writing (see Fig. \ref{fig01}), 
is the integral sign, $D_y$ is the integration domain for $y$, i.e., the $y$-range $[0, 1]$ in Fig. \ref{fig02}. 

If we treat the staircases in Fig. \ref{fig02} as a result of cutting the $x$ domain into small pieces, 
then integration is a sum of vertical small pieces (i.e., the rises of the staircases)  resulting from cutting. 
Note that the cutting, or partition, is not unique. You may make your cutting in your own way. However, to support the 
meaning of ``small pieces" in both $x$- and $y$- domains, the commonly used cuttings or partitions in applications are not too far away from being uniform. 

What is the meaning of $dy$ in the integral in Eq. (\ref{eq0304})? This will be discussed in the next section. We will show that 
\be
dy = f'(x) dx
\ee
where $f'(x)$ is called the derivative function of $f(x)$ and $dx$, like $dy$, symbolizes a cut or partition in the $x$ domain. 
Thus, we have 
\be
 \int_a^b f'(x) dx =  \int_{D_y} \mbox{d}y =  f(b) - f(a).
\ee
This is the fundamental theorem of calculus (FTC). Rather than treating FTC as a theorem in traditional calculus,
we regard it as an axiom or a definition of an integral, an intuitive conclusion based on common sense. 

\section{Meaning of $\Delta y, dy, \Delta x, dx, f'(x)$ and $dy/dx~~~$}
 
\begin{figure} [h]
 \centering
\includegraphics[width=3.5 in]{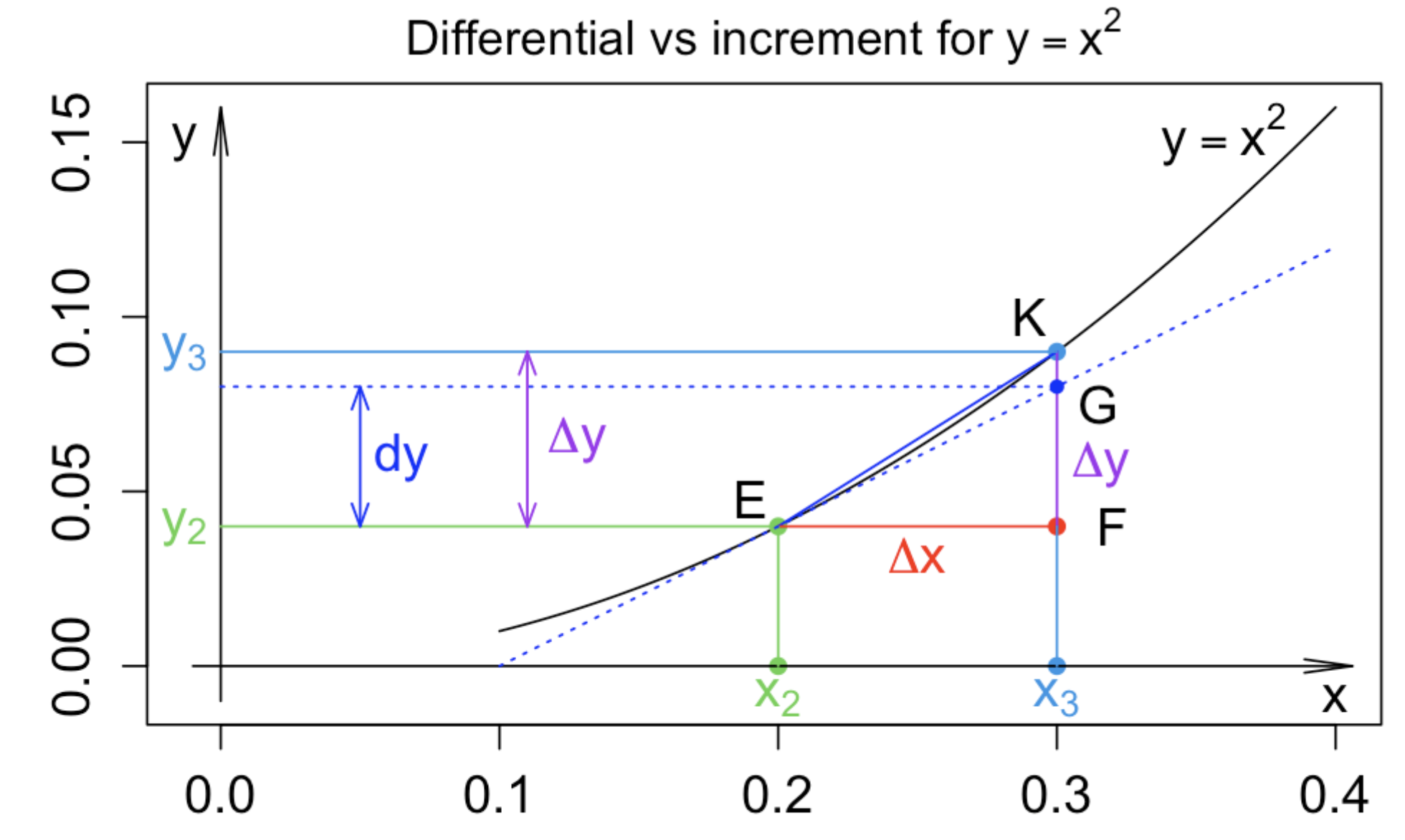}
\caption{Approximate a curve-edged triangle $EFK$ by a right triangle $EFG$.}
\label{fig03}
\end{figure}

 Figure \ref{fig03} is a zoom-in display of Fig. \ref{fig02} for a single piece resulting from the cut. 
 Points $E$ and $K$ are on the curve. 
 For right triangle $\Delta EFK$, the length $FK$ is $\Delta y$, while the length $EF$ is $\Delta x$. 
 The slope of the line segment $\overline{EK}$ is $\Delta y/\Delta x$. 
 
 Next, we explore the 
 slope of the tangent line at point $E$ and its increment $dy = FG$. 
 Assume the slope of the tangent line is $m$.
 Then, the equation for the line is 
 \be
 y - y_2 = m (x - x_2).
 \label{eq3301}
 \ee
 Substitution of $y = x_2$ and $y_2 = x_2^2$ into this equation yields
 \be 
 x^2 - x_2^2 = m (x - x_2)
 \ee
 or 
  \be 
  (x + x_2 - m)(x - x_2) = 0
 \ee
 which has two solutions $x = x_2 $ or $x = -x_2 + m$. Since $E$ is a tangent point, 
 the two solutions must be a repeated root, i.e.,  the two solutions are the same: $x_2 = -x_2 + m$. This implies
 \be
 m = 2 x_2.
 \ee
 This slope is also called the {\it derivative} of the function $y = x^2$ at $x_2$. 
  
 This method of calculating slope without using the concept of limit is from 
 Descrates' direct calculus, called {\it DD calculus},
 in contrast to Fermat's method of tangents, which is similar to our conventional slope 
 definition using the ratio of increments $\Delta y/\Delta x$ (see Section D.1.4 of [2]). 
 Descartes' procedure of finding the slope of $y = x^2$ can be applied to any point $x$,
 resulting in a slope $2x$, which is called the {\it derivative function} of $x^2$, or simply derivative, denoted by 
 $y'$ or $f'(x)$:
 \be
 f'(x)= 2x.
 \ee
 Here, the word ``derivative" may be regarded as a derived quantity. 
 
 The above gives a procedure to calculate the derivative of a function $f(x)$, and hence solves Newton's Problem \#56,
 when $y$ is distance and $x$ is time. 
 
 Although Descartes' direct method for finding the derivative, based on the definition of slope, is conceptually simple and mathematically beautiful, 
 for non-polynomial functions, it is not as practical as Fermat's approximation method of tangents: $\Delta y/\Delta x$ when $\Delta x$ becomes very small. 
 In this sense, we may regard Descartes' direct method as the exact value, while Fermat's method as a numerical approximation. 
 These two methods are bridged by the limit process. 
 
The increment of the tangent line from point $E$ to $G$ is $dy = FG$, which is 
called the differential of $y$ and is equal to the slope times the horizontal increment $\Delta x$:
\be
dy = f'(x) \times \Delta x.
\label{eq3307}
\ee

Next, we show that $dx = \Delta x$. 
The slope of the function $y = x$ is one since its graph is a straight line at 45$^\circ$ from the $x$-axis on the $xy$-plane.
Thus, Eq. \ref{eq3307} implies 
\be
dy = dx = 1 \times \Delta x = \Delta x.
\label{eq3307}
\ee
Thus, you always have $dx = \Delta x$, but $dy$ may not necessarily be $Delta y$. 
We call $dx$ the differential of $x$ when $x$ is specified as the independent variable. 

Figure \ref{fig03} indicates that the difference between $\Delta y$ and $dy$ is equal to the length $KG$.
The difference may be regarded as the error of approximating the function value $f(x)$ at point $K$
by using a tangent line function. We call this a {\it linear approximation} since the tangent line is the 
graph of a linear function. 

The approximation error can be quantified in the following way
\bea
R &=& \Delta y - dy \nonumber \\
&=& \left[(x_2 + \Delta x)^2 - x_2^2 \right] - 2x_2 dx \nonumber \\
 &=&(\Delta x)^2
 \label{eq3312}
\eea
If $\Delta x = 0.1$, then $R = 0.01$, which is much smaller than $\Delta y = 0.05$. 

In summary, we have 
\be
\Delta y = dy + R = f'(x)dx + R.
\label{eq3314}
\ee
Here, the linear approximation error $R$ usually depends on both $x$ and $dx$, although 
$R = (dx)^2$ depends only on $dx$ in the special case of $f(x) = x^2$. 

The quotient of differentials $dy/dx$ is the slope of the tangent line $\overline{EG}$, hence equal to $f'(x)$, and
 is another notation for the derivative. 
The commonly used derivative notations include $y', f', f'(x), dy/dx, df/dx, df(x)/dx$ and $\dot{y}$. The last one 
was coined by Isaac Newton [1], while $dx$ and $dy$ notations were 
attributed to Gottfried Wilhelm Leibniz [3, 4]. 

\section{Interpret $f'(x_i) dx$ and Formulate Applications for Integrals}
We thus have the following approximation
  \be
  f(x)\bigg|_a^b = \int_a^b f'(x) \mbox{d}x  = \sum_{i=1}^n \Delta y_i = \sum_{i=1}^n (f'(x_i)dx + R_i) \approx \sum_{i=1}^n f'(x_i) dx,
  \label{eq0416}
 \ee
 when 
 \be
 \sum_{i=1}^n R_i 
 \ee
 is omitted. The omission is under the condition of small $dx$ because the errors $R_i$ are small. 
 
 The quantity $f'(x_i) dx$ is the area of a rectangle of height equal to $f'(x_i)$ (if $f'(x_i) >0$) and width equal to $dx$. Thus, the geometric 
 meaning of $ \int_a^b f'(x) \mbox{d}x$ is the area underneath of the derivative function $ f'(x)$'s curve bounded by $x=a$, $x=b$ and the $x$-axis. 
 The integral result, by definition, is $f(b) - f(a)$. Namely, the integeral of derivative function $f'(x)$ inside $[a, b]$ can be computed by the original function
 $f(x)$ at the endpoints. The original function $f(x)$ is also called an {\it antiderivative} of $f'(x)$. So, $(f'(x), f(x))$
 form a {\it derivative-antiderivative pair}, or simply DA pair . For example, $(2x, x^2)$ is a DA pair, and $(1, x)$ is 
 also a DA pair. More details of the DA pair concept can be found in Section D.1.1. 3 of [2]. 
 Traditional calculus courses spend much time making hand calculations for 
 the derivative or antiderivative when the other is given. 
 
 It is also the approximation of (\ref{eq0416})
 \be
  \int_a^b f'(x) \mbox{d}x \approx \sum_{i=1}^{n}  f'(x_{i-1}) \Delta x_i
 \ee
 that provides numerous application scenarios for calculus, because many practical 
 applications are set up by formulating a function $g(x)$ times $dx$, 
 similar to $f'(x_{i-1}) \Delta x_i$. For example, $g(x)$ can be the speed of a moving point in the unit of kilometers per hour, which
 is a function of time $x$, then $g(x) dx$ 
 is approximately the distance traveled in a small time interval $[x, x+dx]$. The sum of $g(x) dx$ is the total distance 
 traveled from time $a$ to time $b$. When $dx$ is very small, the sum is equivalent to 
 \be
 \int_a^b g(x)\mbox{d}x,
 \ee
 which is equal to $G(b)-G(a)$, where $G(x)$ is an antiderivative of $g(x)$. 
 
 Another example is to calculate the total construction price of a highway section $[a, b]$ by integration. 
 Because of geological and environmental differences, such as bridges and mountains, the 
 price $p$ is a function of distance $x$ from a base point: $p(x)$. The cost for the highway section $[x, x+dx]$ is approximately 
 $p(x)dx$. The total cost for the highway section $[a, b]$ can be formulated as 
 \be
 \int_a^b p(x) \mbox{d}x.
 \ee

\section{The Mean Value Theorem} 
Since $ \int_a^b g(x)\mbox{d}x$ can be regarged as the area $S$ underneath 
the curve of $g(x)>0$ in $[a, b]$, our common sense concludes that there is a rectangle of the same area with base on $[a, b]$ 
and height $h$ such that 
\be
h = g(c) = \frac{S}{b-a} = \frac{\int_a^b g(x)\mbox{d}x}{b-a}
\label{eq45mvt}
\ee
for at least one point $c$ in $[a, b]$. Here, $h$ is considered an average value, or mean value, of $g(x)$ in $[a, b]$. 
Formula (\ref{eq45mvt}) is conventionally known as the mean value theorem (MVT). 

\begin{figure}[htbp]
    \centering

    \begin{subfigure}[t]{0.4\textwidth}
        \centering
        \includegraphics[width=\linewidth]{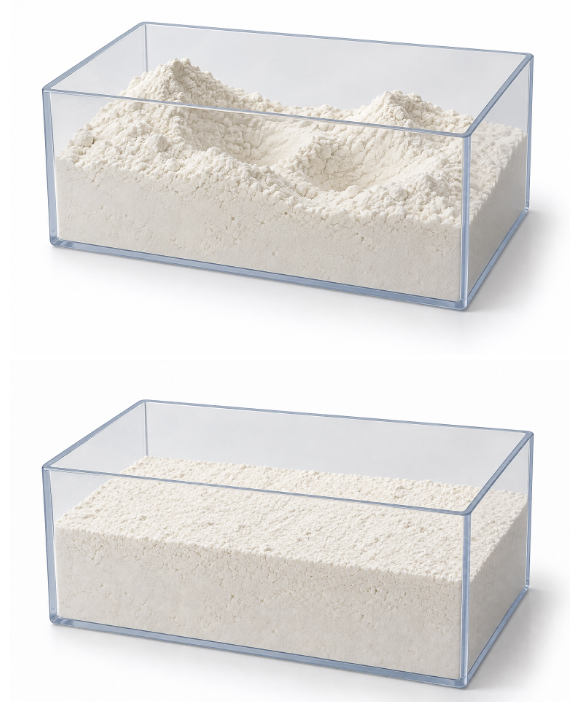}
        \caption{}
    \end{subfigure}
    \hfill
    \begin{subfigure}[t]{0.58\textwidth}
        \centering
        \includegraphics[width=\linewidth]{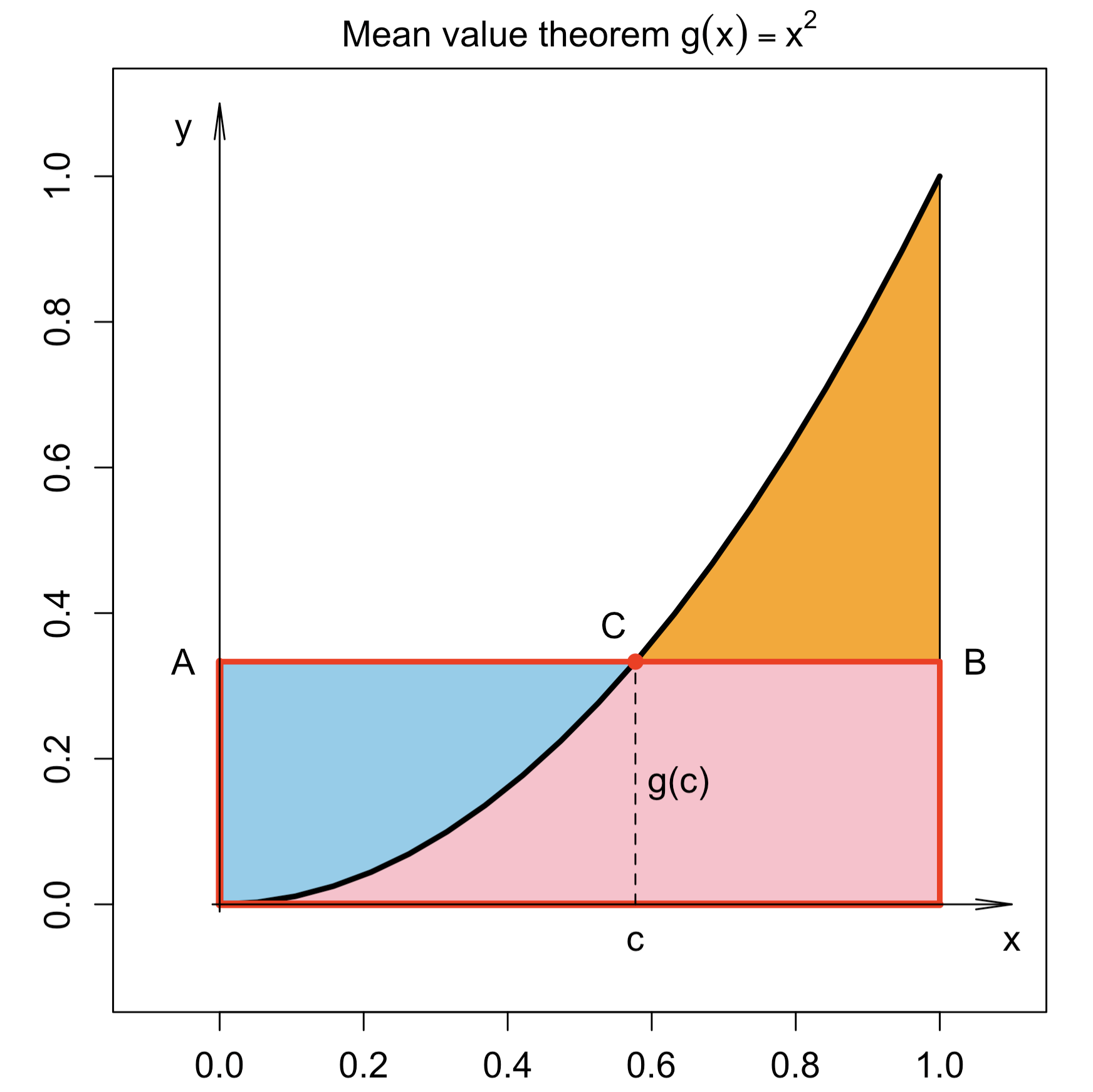}
        \caption{}
    \end{subfigure}

\caption{Graphic demonstration of the mean value theorem (MVT): (a)~ After the food flour is poured into a box container, its surface is heaped. Repeated tapping or vibration redistributes the material and makes the surface nearly level. This leveled surface height is the mean value of the heaped surface.  (b)~The blue deficit area can be filled by the orange surplus area.  The height of the red horizontal line $\overline{AB}$ is the mean value for the function $y = x^2$ in the interval $[0, 1]$.}
\label{fig04}
\end{figure}

MVT is intuitive.  Figure \ref{fig04}(a) shows that an uneven, heaped surface can be leveled by tapping the flour container. 
Figure \ref{fig04}(b) shows that the uneven curve can be leveled by moving the surplus area to the deficit area. 
The area $S$ equals the area of the pink region plus that of the orange region. We can move the horizontal line segment $\overline{AB}$ 
up and down to the height $h$ so that $S = h (b-a)$, i.e., the red rectangle's area is equal to $S$. When this occurs, the blue deficit region can be exactly filled by the orange surplus region. Intuition supports the existence of $h$, and hence the existence of $c$ in $(a, b)$ such that $h = g(c)$. This intuition is the integral form of MVT in today's calculus books. 

The integral MVT can be rigorously proved using the Intermediate Value Theorem (IVT) and Extreme Value Theorem (EVT), both of which
 are based on the completeness axiom of the real numbers. The first time calculus learners do not need to spend time on these proofs. Their concepts and ideas are supported by common sense. 
 IVT says that if Harry wants to walk and meet  Sally on the other side of a river, he must cross the river. This is common sense, and may be regarded as an axiom. MVT is also based on common sense and intuition. It can be regarded as an axiom too.  When assuming one (either IVT or MVT), 
 the other can be proved. 
 
 In calculus teaching, should we treat the MVT as an axiom: the mean value axiom (MVA)? This question may be worth discussion. 
 
 MVA has many important consequences. Two are shown as follows. The first is the differential mean value theorem: 
 \be
 g(c) = \frac{\int_a^b g(x)\mbox{d}x}{b-a} = \frac{G(b) - G(a)}{b-a},
\label{eq46mvtD}
\ee
where $(g(x), G(x))$ is a DA pair. Thus, this formula can also be written as
 \be
  \frac{\int_a^b g(x)\mbox{d}x}{b-a} = \frac{G(b) - G(a)}{b-a} = G'(c).
\label{eq47mvtD}
\ee
This is the differential mean value theorem. Geometrically, this means that a 
secant line between points $(a, G(a))$ and $(b, G(b))$ has at least one parallel tangent line
at point $(c, G(c))$ where $c$ is in $(a, b)$. 

The second is another way to approximate the derivative by a quotient of differentials. In the differential MVT (\ref{eq47mvtD}),
let $a =x$ and $b= x + \Delta x$. MVT becomes
\be
G'(c) = \frac{G(x+\Delta x) -G(x)}{(x+\Delta x) - x} = \frac{\Delta G}{\Delta x},
\ee
where $c$ is in the interval $(x, x + \Delta x)$. This is a more general case of the traditional calculus definition of the derivative 
\be
G'(x) = \lim_{\Delta x \to 0} \frac{\Delta G}{\Delta x}
\ee
when $\Delta x \to 0$ from both right and left sides. 
The MVT's derivative consequence $G'(c)$ is exact, not an approximation. 
Approximations get involved when using  $G'(c)$ to approximate $G'(x)$ under the condition that $c$ and $x$ are close to each other. 
This intuition is easy for students to accept. 

Up to this point, we have solutions to Newton's Problems 56 and 57. For Problem 57, the solution is the FTC:
\be
\int_a^b f'(x) dx = f(b) - f(a).
\label{eq26FTC}
\ee
This is intuitive: $x$ is time, $f(x)$ is speedometer reading at time $x$, $f'(x)$ is speed reading, $ f'(x) dx$ is the 
distance traveled in the time interval $[x, x+ dx]$, and $f(b) - f(a)$ is the difference between the end odometer reading and the start odometer reading.
Because this is intuitive, we suggest using this formulation as the definition of an integral, rather than the area definition or the Riemann sum definition in most textbooks. 

As discussed earlier, the interpretation of $dx, \Delta x, \Delta y, f'(x)$ and $f(x)$ leads to numerous applications (e.g., the area and project cost computations) and numerical computations (e.g., approximating the derivative by a quotient of differentials). The derivative computation solves Newton's Problem 56. 
However, in Ref. [1], Newton did not solve Problem 56 this way; rather, he used the binomial series approach:
\be
(1 + x)^n = 1 + nx + \frac{n(n-1)}{2!}x^2 + \cdots
\ee
For example, he found the derivative of $x^3$ by multiplying its ``arithmetical progression" (i.e., 3), and further by $1/x$. In this way, he solved his Problem 56. 

To prepare for the computation of ``arithmetical progression," his book began with infinite series and devoted the first 18 pages to this subject.  The infinite series was immediately followed by ``Transition to the Method of Fluxions,'' where Problems 56 and 57 were set in italic text (Page 19 of Ref. [1]). The use of italics reflected his emphasis on these core problems for calculus. 
 
 Leibniz, on the other hand, used the approach of differentials and their quotients (See Refs. [3, 4]) to compute slopes. Leibniz's method for computing slopes follows Fermat's method of tangents, while Newton's method follows Descartes' method of tangents. The slopes and tangents of Newton and Descartes are exact, while those of Fermat and Leibniz are approximations. The two approaches are unified by MVT with an exact formulation or by a limit process with approximations. 

\section{Green's Theorem}
Another interpretation of the FTC formula (\ref{eq26FTC}) is prediction. The right-hand side is the antiderivative data at the boundary points of the 
domain $[a, b]$: $f(a)$ and $f(b)$. The integrand under the integral sign on the left-hand side of (\ref{eq26FTC}) is the derivative data inside the domain $[a, b]$: $f'(x)$. 
Thus, one can use the boundary data to predict the integrated 
information inside the domain, or use the inside-domain data to predict the boundary information. Or in general, use surface data to predict the interior information, or use interior data to predict the boundary information. This idea can be applied not only to an interval on the $x$-axis, but also to a 
1D line segment, a 2D region, a 3D body, or even a higher-dimensional bounded domain. 

For an interval on the $x$-axis, this prediction may be considered a special case of Newton's Problem 56 with less required data. 
Instead of ``The length of the Space described being continually (that is, at all
Times) given," this prediction formulation gives only the length at two times: $x=a$ and $x=b$. What is lost is the uniqueness of the solution to Problem 56: 
the speed inside the interval $[a, b]$ can have infinitely many possibilities. This is like when driving from San Diego to Los Angeles from 2 pm to 4 pm; each driver may have different speeds at different times, although they all start at 2 pm and arrive at 4 pm.  MVT provides only the average speed
\be
v = \frac{f(b) - f(a)}{b - a}
\ee
and every driver must at one point have driven at this average speed $v$.

The condition for this prediction problem is still the same as that of Newton's Problem 57: $f'(x)$ is given for every point inside $[a, b]$. 
Then, the difference of the boundary information $B_f = f(b) - f(a)$ can be predicted by an analytic method or by numerical integration 
\be
B_f =  \int_a^b f'(x) dx \approx \sum_{i=1}^n f'(x_i) \frac{b-a}{n}, 
\ee
where $n$ is the number of uniform partitions of $[a, b]$. 

\begin{figure} [h]
\centering
\includegraphics[width=3.1in]{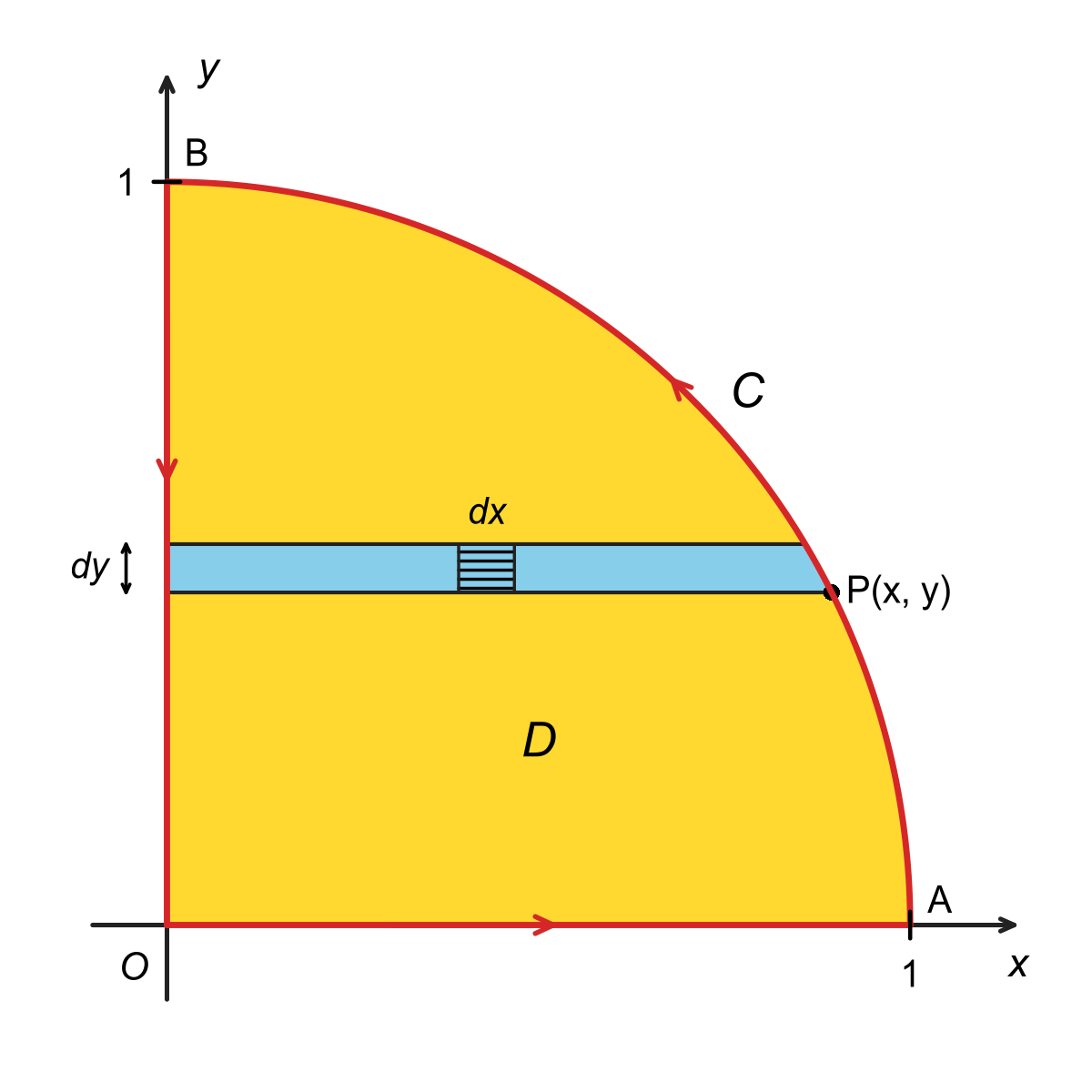}
\caption{A quarter of the unit disc is used to support the explanation of Green's theorem. The small partitioned area is $dA = dx \times dy$ in this case. We first cut a 
horizontal strip of width $dy$ out of the domain $D$. Then cut the strip vertically with a width $dx$. The point $P$ with coordinates $(x, y)$ 
satisfy $x = \sqrt{1 - y^2}$ because of the unit circle.}
\label{fig05}
\end{figure}

This prediction formulation corresponds to many natural phenomena. For example, moving a magnet through a region bounded by a closed conducting wire 
generates electricity in the wire. This is the principle of electrical power generation in our power stations. The electric flow can be predicted by the magnet's strength and its speed. More precisely, the electric current can be predicted by the flux of the changing magnetic field. Mathematically, the inside flux is the derivative function, and the boundary is the original function:
\be
\int_D f_x ~dA = \oint_C f ~dy,
\label{eq71Green1}
\ee
where $f$ is a function of two variables $(x, y)$ in a 2D region $D$ (See Fig. \ref{fig05}), $dA$ is a small area resulting from the partitition of $D$,
$f_x$ is the derivative function with respect to $x$ (or called the x-partial derivative), $dy$ is diffential on the curve $C$, which is the boundary of the region $D$ and has counterclockwise direction. The notation $\oint$ indicates that the integration domain $C$ is a closed curve. 
The curve $C$ is closed because it encloses a region and conducts electricity. It is also natural for the curve $C$ to have a direction, as electric flow has direction. The counterclockwise direction corresponds to Lenz's law in physics, which determines the direction of the electric current. 

The integral over a surface or a region is called a surface integral, like the left-hand side integral in Eq. (\ref{eq71Green1}), and that over a curve is called 
line integral or curve integral, like the right-hand side integral in Eq. (\ref{eq71Green1}). 

Equation (\ref{eq71Green1}) is known as Green's theorem, attributed to British mathematical physicist George Green (1793 - 1841), who published the theory in 1828 in his research on 
mathematical analysis for electricity and magnetism. This equation shows that 
the amount of the derivative function going through a 2D region can predict the circulation of the original function on the boundary of the region. 
In Green's theorem, both sides are expressed as integrals. For FTC, we can also regard $f(b) - f(a) = f\big|_a^b$ as an integral over the
0D domain consisting of two points at $a$ and $b$. With this, we see that Green's theorem reduces a derivative function's integral over a 2D domain to an original function's integral over a 1D domain, while the FTC reduces a derivative function's integral over a 1D domain to an original function's integral over a 0D domain. Naturally, you might ask: can this be extended to 3D, 4D, or even higher-dimensional space? The answer is yes. The 3D case is the divergence theorem to be discussed in the next section. 

Next, we prove Green's theorem (\ref{eq71Green1})  based on the 2D domain shown in Fig. \ref{fig05}. 

Apply the FTC to the left-hand side of Eq. (\ref{eq71Green1}) on the blue horizontal strip for $f_x dx$ with $dA = dx\times dy$
\bea
\int_D f_x ~dA = \int_D f_x ~dx\times dy &=& \int_{D_y} f(x, y)\bigg|_0^{\sqrt{1-y^2}} ~dy \nonumber \\
&=& \int_{D_y} \left[f(\sqrt{1-y^2}, y)  - f(0, y)\right] ~dy,
\label{eq72}
\eea
where $D_y = [0, 1]$ is the integration domain (or called integration interval) for $dy$. 

The integration
\be
\int_{D_y} f(\sqrt{1-y^2}, y)~ dy = \int_0^1 f(\sqrt{1-y^2}, y)~ dy =  \int_{\overset{\frown}{AB}} f(x, y) ~dy
\ee 
is along the arc $\overset{\frown}{AB}$ from $A$ to $B$ following the direction indicated by the red arrow on the arc.

The integration
\be
-\int_{D_y} f(0, y)~ dy = -\int_0^1 f(0, y)~ dy = \int_1^0 f(0, y)~ dy = \int_{\overline{BO}} f(x, y) ~ dy
\ee 
is along the line segment $\overline{BO}$ from point $B$ to point $O$, also following the direction indicated by the red arrow on this line segment. 

Along the line segment $\overline{OA}$, $dy = 0$ because a horizontal line cannot have a $y$-increment. Thus, 
\be
\int_{\overline{OA}} f(x, y)~ dy = 0.
\ee 

Substitution of the above three integrals into Eq. (\ref{eq72}) yields
\be
\int_D f_x ~dA  =  \int_{\overset{\frown}{AB}} f(x, y) ~dy + \int_{\overline{BO}} f(x, y) ~ dy + \int_{\overline{OA}} f(x, y)~ dy. 
\ee
The integration domains $\overset{\frown}{AB}, \overline{BO}$ and $\overline{OA}$ of the three integrals on the right-hand side form a closed curve $C$, 
 which is in 
counterclockwise direction and leaves its bounded domain to the left. Therefore, the sum of these three integrals can be 
written as 
\be
\oint_C f(x,y) ~dy, 
\ee
the right-hand side of Green's theorem (\ref{eq71Green1}). This completes the proof. 

Although this proof of Green's theorem used a specific domain: a quarter-disc. The proof is valid for any 2D domain, as a connected 2D domain can be partitioned into several domains that look like a quarter-disc, although the arc side may not be round. The first-time learners of calculus, of course, do not 
need to spend too much time figuring out the details; instead, you may focus on the basic ideas. When you review the subject, you can explore more details with more examples and real-world applications. Our goal here is to help instructors teach calculus in such a way that students never forget the basic ideas like FTC, MVT, and Green's theorem. They may not remember the names, but they will not forget these common-sense-based ideas. 

As an application, Green's theorem can compute the area of a lake from coordinate data on the lake's boundary. Let $f(x,y) = x$, then $f_x = 1$. Thus,
the left-hand side of Green's theorem (\ref{eq71Green1}) is
\be
\int_D f_x ~dA = \int_D 1 ~dA = A,
\ee
which is the area of the 2D region $D$ and is equal to the right-hand side 
\be
A = \oint_C x ~dy \approx \sum_{i=1}^n x_i \Delta y_i.
\ee
Here, $x_i$ the $i$th point's $x$-coodinate and $\Delta y_i$ is the $i$th point's $y$-increment, after cutting the boundary curve 
$C$ into $n$ pieces. As long as the boundary coordinate data $(x_i, y_i)$ are given, the area $A$ can be computed. 
You do not need to ride a boat to measure the lake surface directly!

This formula is intuitively true as $x\times dy$ is the area of the horizontal blue strip in Fig. \ref{fig05}. Thus, the entire area is an integral of $x\times dy$.

The planimeter apps for computing an area bounded by a closed curve can be downloaded from the App Store 
for iPhones or the Google Play Store for Android phones and are programmed based on Green's theorem. 
Planimeters have important medical uses, such as computing the area of the cross-section of a bone or an organ for 
operation or diagnosis without cutting or damaging the bone or the organ. 

You may ask why the partial derivative in Green's theorem (\ref{eq71Green1}) is with respect to $x$? Does it work for $f_y$?
Yes, it should. And it does. The corresponding formula is 
\be
\int_D f_y ~dA = - \oint_C f ~dx.
\label{eq71Green2}
\ee
The proof is similar to that for $f_x$. The negative sign here comes from the counterclockwise direction of the curve $C$.

A critically important application of Green's theorem is to vector fields. In a 2D space, every point may have a vector,
such as a wind vector field over a lake surface. The vector field may be denoted by ${\bf u}({\bf r})$, where ${\bf u} = (u, v)$
and ${\bf r} = (x, y)$ is called the position vector. So a vector field is a function defined at position vectors. 

Some vector fields rotate, such as the wind field in tornadoes and hurricanes. The strength of the rotation is measured by vorticity, defined by 
\be
\omega = v_x - u_y. 
\ee
Based on Eqs. (\ref{eq71Green1}) and (\ref{eq71Green2}), Green's theorem for vorticity $v_x - u_y$ is as follows
\be
\int_D \left(v_x - u_y\right)  ~dA =  \oint_C u~dx + v~ dy = \oint_C {\bf u} \cdot d{\bf r}
\label{eq71Green3}
\ee
where the notation $\cdot$ is for a dot product of two vectors defined as follows:
\be
{\bf u} \cdot d{\bf r} = u~dx + v~ dy.
\ee
Equation (\ref{eq71Green3}) is the usual presentation of Green's theorem in most of the existing textbooks. 
This formulation is difficult for many students. Most students cannot remember this formula shortly after the final exam. 

On the other hand, our approach presents Green's theorem as an extension of the FTC from 1D to 2D and should be easier for students to 
comprehend. Expressions (\ref{eq71Green1}) and (\ref{eq71Green2}) are also simpler than Eq. (\ref{eq71Green3}). 

In electromagnetism, if ${\bf E} = {\bf u} = (u, v)$ is the electric field, then its vorticity is the rate of change of the corresponding 
magnetic field. The mathematical expression is 
\be
\int_D B_t ~dA = \oint_C {\bf E}\cdot d{\bf r},
\label{faraday1}
\ee
where $B$ is the magnetic field's component in the direction perpendicular to the region $D$ whose boundary is the closed curve $C$, 
 $B_t$ is the rate of change of the magnetic field,  ${\bf E} = (E1, E_2, 0)$ is the electric field along $C$, and 
 \be
 B_t = \frac{\partial E_1}{\partial y} - \frac{\partial E_2}{\partial x}.
 \label{faraday2}
 \ee
Equation (\ref{faraday1}) is the integral form of Faraday's law on a flat surface $D$, one of the four Maxwell equations. 
Equation (\ref{faraday2}) is the differential form of Faraday's law. The integral form and the differential form are related by 
Green's theorem with a negative sign, which is due to Lenz's law on the direction of the electric field ${\bf E}$ induced by the changing 
magnetic flux through $D$.

\section{The Divergence Theorem}
\begin{figure} [h]
\centering
\includegraphics[width=5.0in]{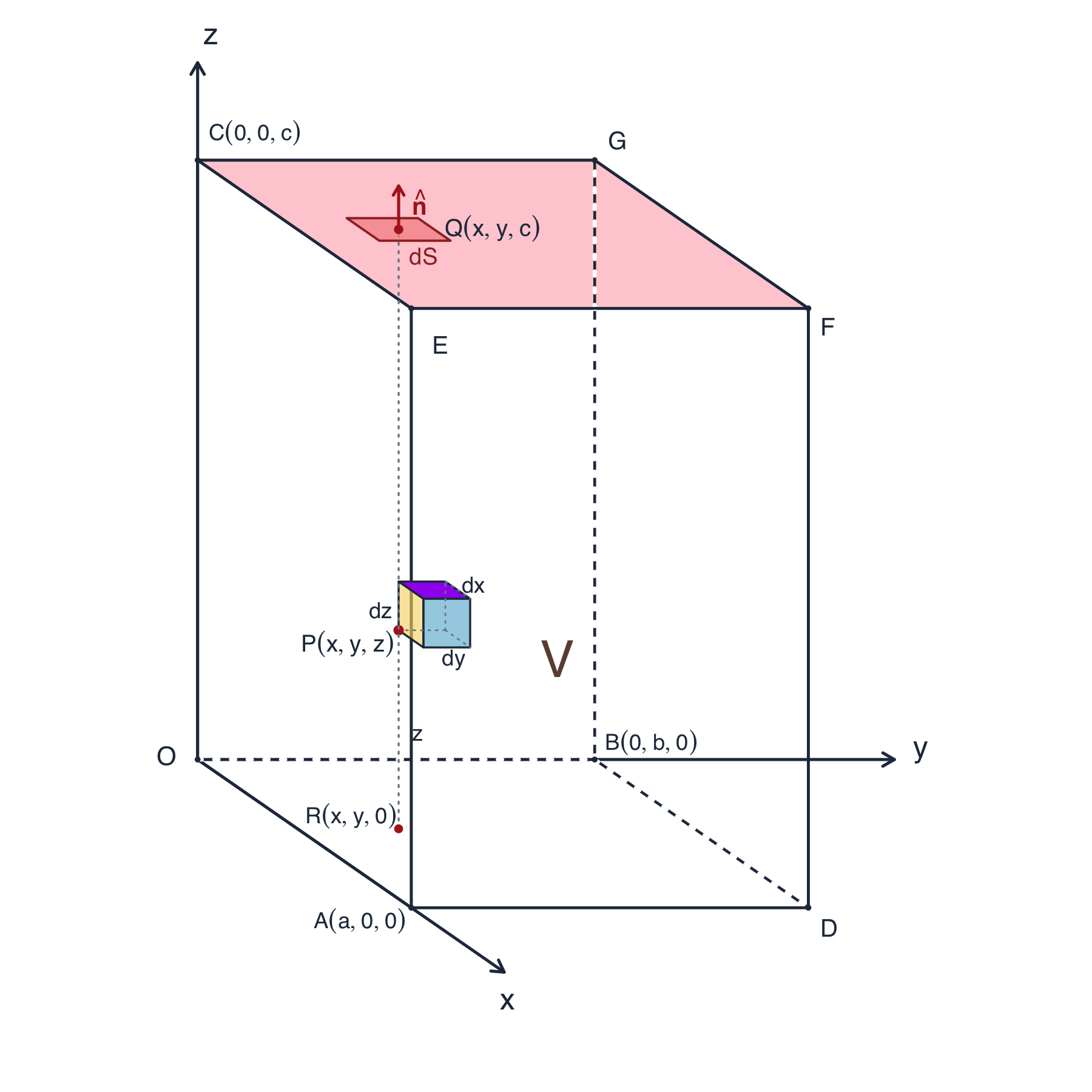}
\caption{This vertical rectangular block, or called a prism, is for the proof of the divergence theorem. The small volume is $dV = dx\times dy \times dz$}
\label{fig06}
\end{figure}

FTC's extension to a 3D domain is the divergence theorem: The integration of a derivative function inside a domain $V$ 
can be predicted by integrating the corresponding original function on the boundary, denoted by $\partial V$, i.e., 
\be
\int_V f_z ~dV = \oint_{\partial V} f \times n_3 ~dS,
\label{eq44diver1}
\ee
where $f(x, y, z)$ is a function defined on a 3D domain $V$, $f_x$ is its $x$-partial derivative, $dS$ is a small area resulting from a 
partition of the closed surface denoted by ${\partial V}$,  and the unit 
normal vector for the small area $dS$ is denoted by ${\bf n} = (n_1, n_2, n_3)$ where $n_1, n_2$ and $n_3$ are the $x-$, $y-$ and $z-$ components of the vector ${\bf n}$.
See Fig. \ref{fig06}. 

We prove the divergence theorem (\ref{eq44diver1}) for a specific 3D domain $V$: a vertical rectangular block. 
Similar to the idea of Green's theorem for a general 2D domain, a general 3D domain may be partitioned into 
several blocks, and the divergence theorem is applicable to each block. 

Applying FTC to $dz$ in the left-hand side of the divergence theorem (\ref{eq44diver1}) yields
\bea
\int_V f_z ~dV &=& \int_V f_z ~dx\times dy \times dz \nonumber \\
&=& \int_V \left(f_z ~  dz\right) dx\times dy \nonumber \\
&=& \int_D \left(f_z(x, y, c) - f_z(x, y, 0)\right) dx\times dy \nonumber \\
&=& \int_{D_{top}} f_z(x, y, c)~ dS - \int_{D_{bottom}}f_z(x, y, 0) ~ dS \nonumber \\
&=& \int_{D_{top}} f_z(x, y, c)~ n_3 dS + \int_{D_{bottom}}f_z(x, y, 0) ~ n_3 dS,
\label{eq45div}
\eea
since the unit normal vector ${\bf n}$ are $(0,0,1)$ for the top face $D_{top}$, and  $(0,0,-1)$ for the bottom face $D_{bottom}$.

On the four side faces of $\partial V$, $n_3 = 0$, hence $n_3 dS = 0$ and 
\be
\int_{D_{sides}}f_z(x, y, 0) ~ n_3 dS = 0
\ee

Thus, the expression (\ref{eq45div}) can be written as 
\be
\int_{D_{top}} f_z(x, y, c)~ n_3 dS + \int_{D_{bottom}}f_z(x, y, 0) ~ n_3 dS + \int_{D_{sides}}f_z(x, y, 0) ~ n_3 dS.
\ee
 This sum of integrals is equal to 
 \be
 \oint_{\partial V} f \times n_3 ~dS,
 \ee
because the surface integration domains include all the faces of the rectangular block, i.e., 
\be
D_{top}\cup D_{bottom}\cup D_{sides} = \partial V.
\ee
The proof is complete. 

Note that we purposely used $dS$ to indicate a small surface area on $\partial V$, instead of $dx\times dy$, because the surface 
integral $\oint_{\partial V}$ can have its own choices of the domain partition of the surface $\partial V$, different from the partition of the 3D domain $V$. 

Similar to Eq. (\ref{eq44diver1}), the expressions for the divergence theorem with respect to $x$ and $y$ derivatives can be written as
\bea
\int_V f_x ~dV &=& \oint_{\partial V} f \times n_1 ~dS, \label{eq44diver2}\\
\int_V f_y ~dV &=& \oint_{\partial V} f \times n_2 ~dS. \label{eq44diver3}
\eea

When applying the divergence theorem to a vector field in 3D, such as the water flow current ${\bf u} = (u_1, u_2, u_3)$, using the 
three divergence theorem expressions (\ref{eq44diver1}), (\ref{eq44diver2}), and (\ref{eq44diver3}) for the $x-, y-$ and $z-$ derivatives, we have
\be
\int_V \left[ (u_1)_{x}  + (u_2)_{y} + (u_3)_{z} \right] ~dV = \oint_{\partial V} \left( u_1  n_1 + u_2  n_2 + u_3  n_3\right) ~dS.
\label{divthm}
\ee
This is the divergence theorem presented in most textbooks. We choose to present the theorem using a single component $f_z$ first, following the OFFER
idea: Everything in calculus is built on one fundamental formula, the FTC. 

The expression $(u_1)_{x}  + (u_2)_{y} + (u_3)_{z} $ has important physical meaning: divergence of the vector field ${\bf u} = (u_1, u_2, u_3)$. 
Opposite to divergence is convergence. Divergence means flowing out like a source of a water fountain, and convergence means flowing in like a water
sink in a bathtub. See Fig. 7(a) for a diverging vector field, and (b) for a converging vector field in a 2D space. 

 \begin{figure} [h]
\centering
\includegraphics[width=2.2in]{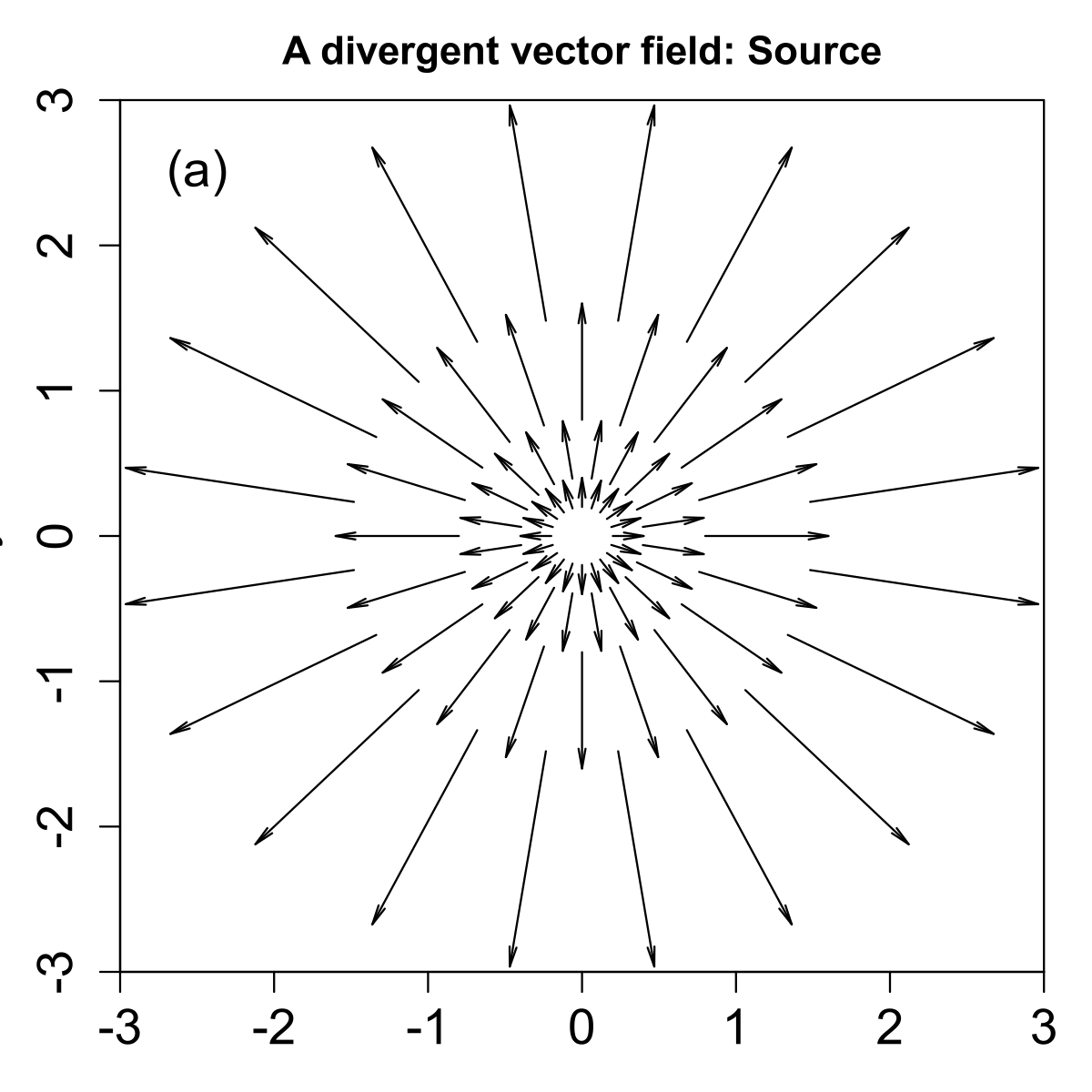}
\includegraphics[width=2.2in]{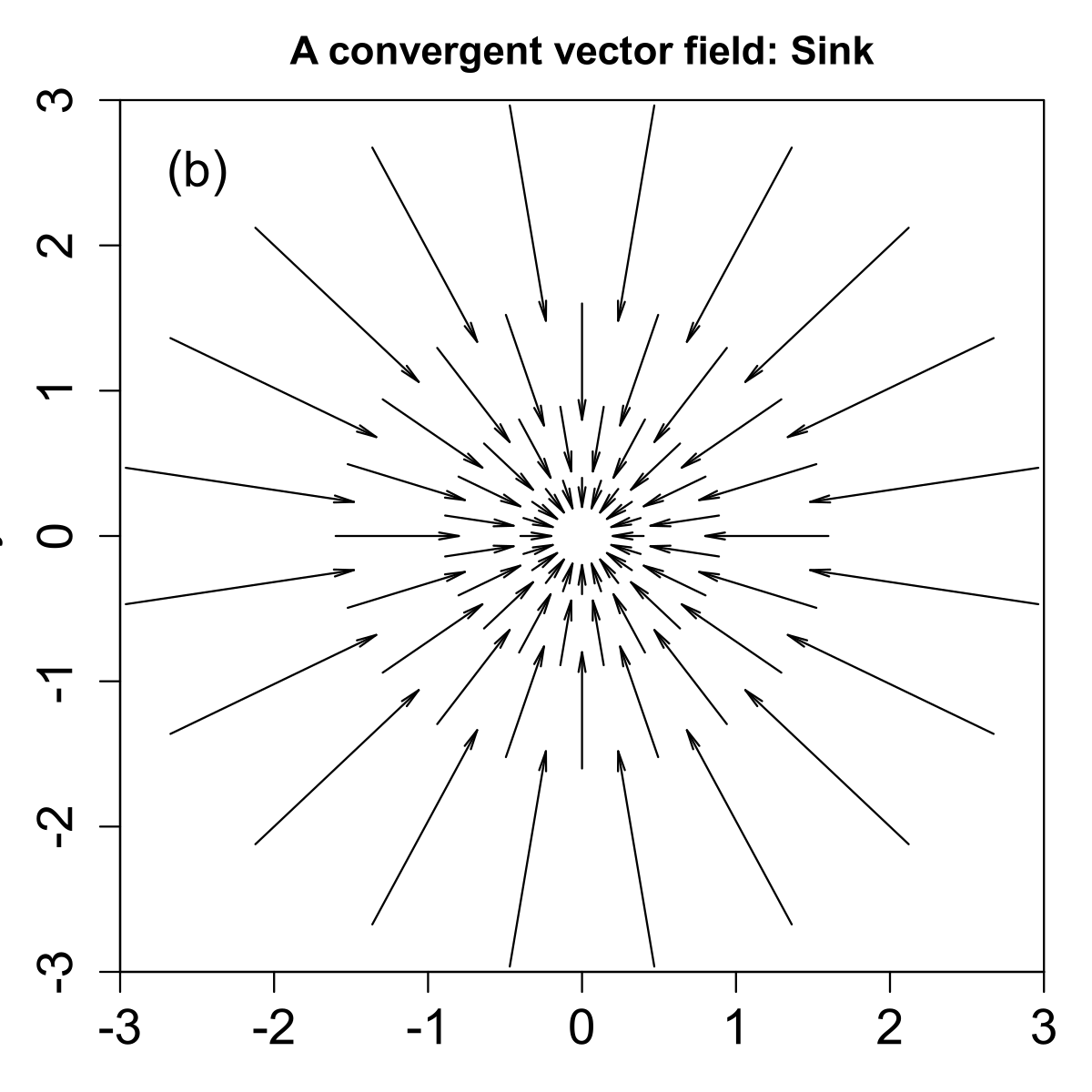}
\caption{(a) ~ A divergent vector field: ${\bf u} = (x/2, y/2, 0)$. (b)~ A convergent vector field: ${\bf u} = (-x/2, -y/2, 0)$.}
\label{fig07}
\end{figure}

Physics often uses a vector notation for divergence 
\be
\mbox{div} ~ {\bf u} = \nabla \cdot {\bf u},
\ee
where the notation
\be
\nabla = \left( \frac{\partial }{\partial x},  \frac{\partial }{\partial y},  \frac{\partial }{\partial z} \right).
\ee
is pronounced as {\tt nabla} or {\tt del}. 

When $\mbox{div}~ {\bf u} ({\bf r}) > 0$, the vector field ${\bf u} $ is divergent at the location ${\bf r} = (x, y, z)$. 
When $\mbox{div}~ {\bf u} ({\bf r}) < 0$, the vector field ${\bf u} $ is convergent at the location ${\bf r} = (x, y, z)$. 
When $\mbox{div}~ {\bf u} ({\bf r}) =0$ everywhere, the vector field ${\bf u} $ is said sourceless. 

In electromagnetism, the electric field $ {\bf E}$'s divergence in a 3D domain is related to the electric charge inside the domain:
Positive divergence corresponds to positive charge, and vice versa. The divergence theorem for a charged 3D domain $V$ and its surface $\partial V$
is Gauss's law of electricity, one of the four Maxwell equations:
\be
 \int_V \frac{q}{\varepsilon_0} ~dV= \oint_{\partial V} {\bf E} \cdot {\bf n} ~dS,
\ee
where 
\be
q({\bf r}) = \nabla \cdot {\bf E}
\ee
 is the electric charge density, $\varepsilon_0$ is the electric constant, also called the vacuum permittivity,
 and $\partial V$ is called a Gaussian surface in physics. 

\section{A General Domain-Boundary Prediction Problem in 3D Space}

\begin{figure} [h]
\centering
\includegraphics[width=4.5in]{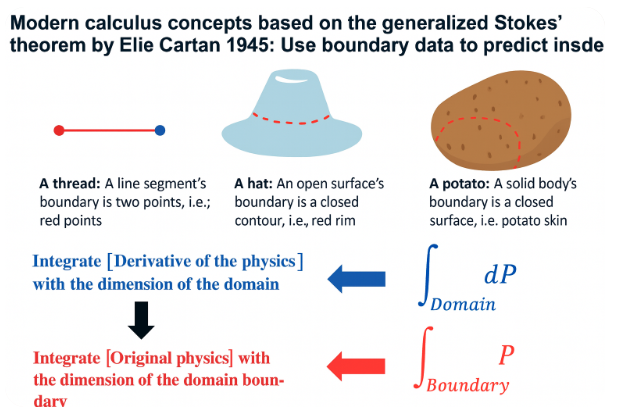}
\caption{The interior and boundary of the three kinds of geometric objects in 3D space: a curve, a surface, and a solid body. 
This version of the figure was generated with the AI tool ChatGPT from the author's original hand-drawn version.}
\label{fig08}
\end{figure}

In 3D space, there are only three kinds of geometric objects: a curve measured by length, a surface measured by area, and 
a body measured by volume. See Fig. \ref{fig08}. 

A curve $L$, measured by length, is a 1-dimensional object since 
it can be determined by three parametric equations with a single parameter, which can be time $t$ or the length $\ell$ of the curve from an end. The boundary of the curve $L$ consists of the two endpoints, which are 0-dimensional. 

An open surface $S$, measured by area, is a 2-dimensional object since it can be determined by two parameters, such as latitude and longitude. Its boundary $\partial S$ is a closed curve, such as the rim of a hat or the bank of a lake. 

A solid body $V$, measured by volume, is a 3-dimensional object, as its determination requires three parameters, such as latitude, longitude, and altitude. 
Its boundary $\partial V$ is a closed surface, such as the skin of a potato or the surface of the Earth. 

The FTC prediction idea applies to any of these three geometric objects. The applications have profound mathematical, physical, and engineering consequences, such as the theory of gravitational potential and Maxwell's equations. 

This interior-boundary prediction can be expressed as
\be
\int_D \delta f  = \oint_{\partial D} f 
\label{eq81sc}
\ee
where $\delta$ denotes specific derivatives of the function $f({\bf r})$, depending on whether the function $f$ is a 
scalar or a vector and on whether the interior domain is a curve, a surface, or a solid body. 

\subsection{Gradient of a Function and its Integration along a Curve}
For a scalar function $\phi$ on a curve $L$ determined by
parametric equation ${\bf r} = {\bf r}(t)$ starting at point $A$ and ending at $B$, $\delta f$ is
\be
\delta f = \nabla \phi \cdot d{\bf r},
\ee
where $\phi({\bf r})$ is called the potential function (e.g., gravitational potential, electrical potential), $\nabla \phi$ is the 
gradient of the potential (also denoted by $\mbox{grad} f$), and $d{\bf r}$ is the differential of the curve $L$.  Then, the prediction equation (\ref{eq81sc}) for the potential
$\phi$ and curve $L$ becomes
\be
\int_L \nabla \phi \cdot d{\bf r} = \phi({\bf r})\bigg|_{{\bf r}_A}^{{\bf r}_B} = \phi({\bf r}_B) - \phi({\bf r}_A),
\ee
where ${\bf r}_A$ is the position vector for point $A$, and ${\bf r}_B$ for point $B$. For the gravitational potential, 
this equation means that the work done by the gravitational force $\nabla \phi $ is $\int_L \nabla \phi \cdot d{\bf r} $ and is equal to 
the difference in the potential. This is consistent with the energy conservation law. 

As discussed earlier, integration along a curve like $\int_L$ is called a line integral. 
Similarly, integration over a surface is called a surface integral, and integration over a 3D solid body is called a volume integral or body integral. 

\subsection{Curl and its Integration on an Open Surface}

When the function $f$ in Eq. (\ref{eq81sc}) is a vector ${\bf u} = (u, v, w)$, there are two cases of $\delta f$. One is the curl on a surface
\be
\delta f = \nabla \times {\bf u} \cdot {\bf n} dS.
\label{eq83curl}
\ee
Another is the divergence over a volume
\be
\delta f = \nabla \cdot {\bf u} ~dV.
\label{eq84div}
\ee

This subsection describes the curl case, and the next describes the divergence. 

Curl, by definition, means rotation or vorticity. Figures 9(a) and (b) show a vector field with positive and negative curl, respectively. 

 \begin{figure} [h]
\centering
\includegraphics[width=2.2in]{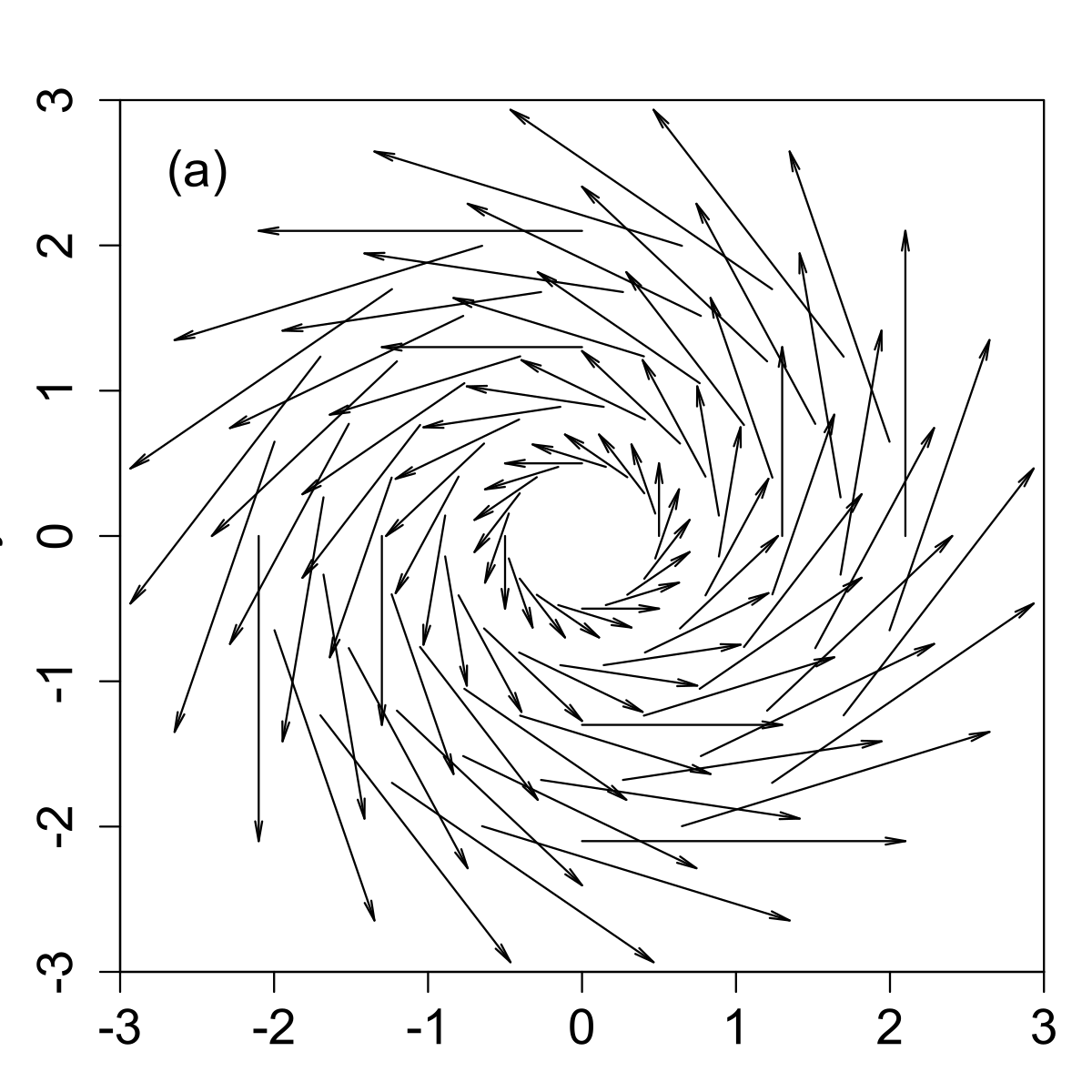}
\includegraphics[width=2.2in]{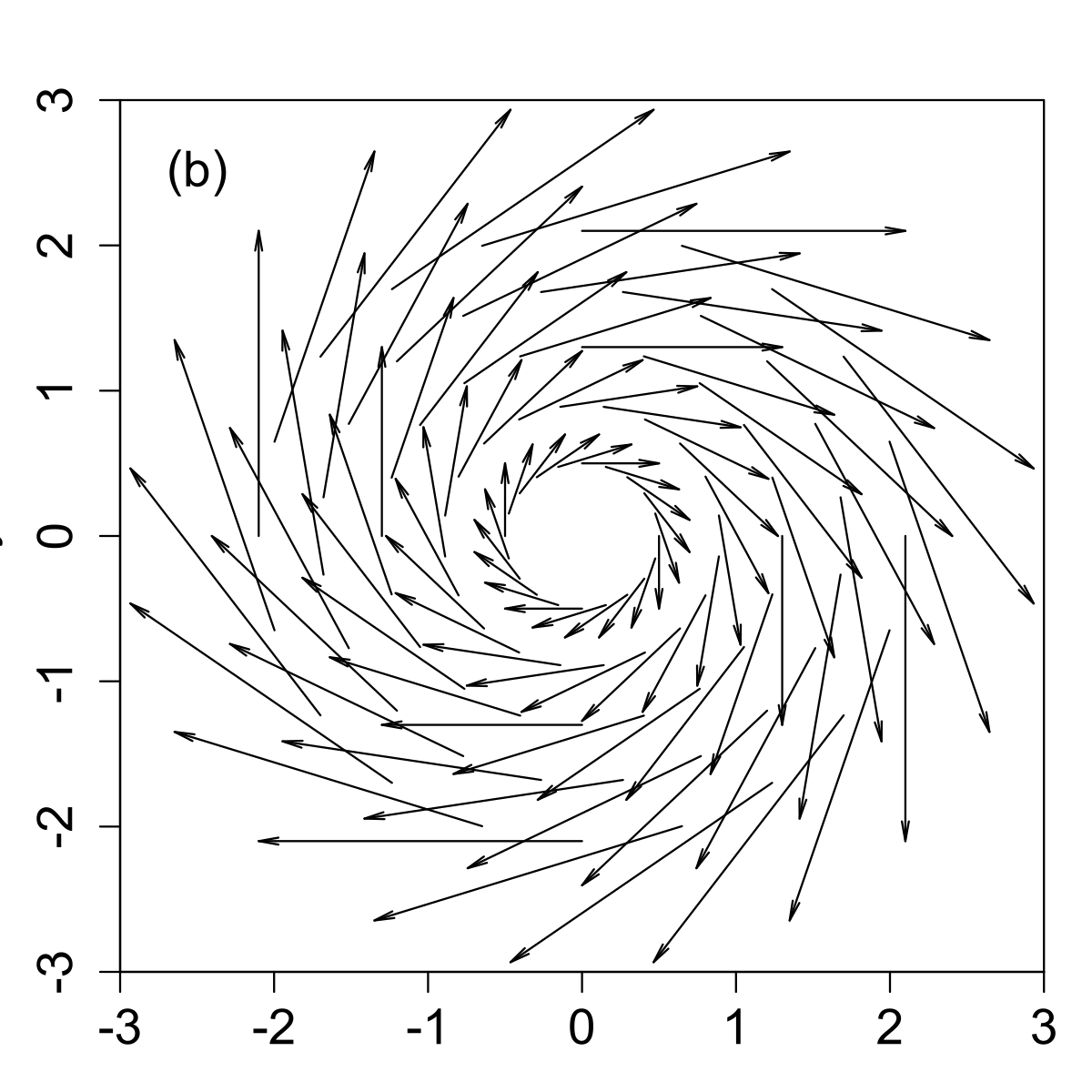}
\caption{(a) ~ A vector field with positive curl rotating counterclockwise: ${\bf u} = (-y, x, 0)$. (b)~ A vector field
with negative curl rotating clockwise: ${\bf u} = (y, -x, 0)$.}
\label{fig09}
\end{figure}

The interior-boundary prediction equation (\ref{eq81sc}) for an open surface $S$ becomes
\be
\int_S \nabla \times {\bf u} \cdot {\bf n} dS = \oint_{\partial S} {\bf u} \cdot d{\bf r}.
\label{eq83curl2}
\ee
Here, 
\be
\mbox{curl} ~ {\bf u}  = \nabla \times {\bf u}
 = \left(
\frac{\partial w}{\partial y} -\frac{\partial v}{\partial z}, 
\frac{\partial u}{\partial z} - \frac{\partial w}{\partial x},
\frac{\partial v}{\partial x} - \frac{\partial u}{\partial y} 
\right)
\ee
The left-hand side of Eq. (\ref{eq83curl2}) means a summation of the local vorticities over $dS$, which results in a
global circulation along the boundary of $S$ given by the right-hand side. Here, the global circulation
is the integral of the vector field along a closed curve, which in this case is $\partial S$. 
 Intuitively, this result also makes sense: the sum of all the local circulations inside a closed curve is equal to the total circulation
 along the curve. 
 
 This equation is Stokes's theorem, attributed to the British mathematician and physicist George Stokes (1819-1903). 
Green's theorem, discovered earlier in 1828, was for a flat surface $S$. When $S$ is lifted to become a heaped surface, Green's theorem becomes Stokes theorem. William Thomson (Lord Kelvin) (1824-1907), a Scottish mathematician,  had this intuition and insight. In his letter to Stokes in 1850, he explicitly
stated the theorem in modern form
\bea
&&\int (\alpha dx + \beta dy + \gamma dz)  \nonumber \\
&=&\pm \iint \left\{ 
\ell \left(\frac{d\beta}{dz} -\frac{d\gamma}{dy} \right)
+  m \left(\frac{d\gamma}{dx} -\frac{d\alpha}{dz} \right)
+ n \left(\frac{d\alpha}{dy} -\frac{d\beta}{dx} \right)
\right\}.
\eea
 George Stokes included it as a problem on the 1854 Cambridge Smith’s Prize examination and asked students to prove it; it later became known as Stokes’ theorem. Rigorous proofs followed a few years later. 
 
 Faraday's law (\ref{faraday1})  also holds for a non-flat surface $S$ and can be expressed as 
 \be
 -\int_S {\bf B}_t \cdot {\bf n} ~dS = \oint_{\partial S} {\bf E}\cdot d{\bf r},
 \label{faraday3}
 \ee
 and 
 \be
 -{\bf B}_t = \nabla \times {\bf E}.
  \label{faraday3d}
 \ee
 Here, the negative sign is due to Lenz's law for the direction of the electric current induced by the 
 changing magnetic flux through the surface $S$. 
 Similar to those discussed earlier on Green's theorem and Faraday's law 
 for a flat region $D$, Eq. (\ref{faraday3}) is the general 3D version of Faraday's 
 law in integral form, and Eq. (\ref{faraday3d}) in differential form, one of the four
 Maxwell equations. 

\subsection{Divergence and its Integration over a 3D Solid Domain}

The interior-boundary prediction equation (\ref{eq81sc}) for a 3D solid body $V$ becomes
\be
\int_V \nabla \cdot {\bf u} ~dV = \oint_{\partial V} {\bf u} \cdot {\bf n} ~dS.
\label{divthmvec}
\ee
This is the vector form of the divergence theorem given by Eq. (\ref{divthm}). Here, the 3D solid body $V$ may be 
regarded as an open domain in 4D space. 

The appearance of the vector form of the divergence theorem (\ref{divthmvec})  is 
 convenient not only in mathematical writing (e.g., a proof of the theorem for a general domain $V$, 
rather than the special rectangular prism shown in Fig. \ref{fig06}), but also in applications in 
physics and engineering (e.g., derivation of the mass conservation equation in fluid dynamics). 
It is also geometrically intuitive. 

In physics and engineering, the divergence theorem for a volume is often related to the law of mass conservation, the FTC along an open curve 
is related to the law of energy conservation, and Stokes' theorem for an open surface is related to the law of vorticity conservation. 


\section{Taylor's Theorem}
A mathematical application of the FTC is deriving Taylor's theorem. 
In calculus education, proving Taylor's theorem can be difficult for many students; so is understanding the reasoning behind Taylor series. 
Here, we show that Taylor's series is a result of the repeated use of the FTC, as shown in Section 6.8.1 of Ref. [2]. 

FTC implies that
\bea
f(x) &=& f(0) + \int_0^x dt_1~ f'(t_1)  \nonumber \\
&=& f(0) +\int_0^xdt_1~ \left[f'(0)+ \int_0^{t_1}dt_2~f''(t_2) \right] \nonumber\\
&=& f(0)+\int_0^xdt_1 f'(0) + \int_0^x dt_1 \int_0^{t_1}dt_2~f''(t_2)  \nonumber\\
&=& f(0)+f'(0)x+ \int_0^x dt_1 \int_0^{t_1}dt_2~f''(t_2). 
\eea
The sum of the first two terms $f(0) + f'(0)x$ is a linear approximation of $f(x)$ in the neighborhood of $0$. 

Continuing the FTC application to the above equation for $f''(t_2)$ , we have 
\bea
f(x) &=& f(0)+f'(0)x+ \int_0^x dt_1 \int_0^{t_1}dt_2~\left[f''(0)+ \int_0^{t_2}dt_3~f'''(t_3) \right] \nonumber \\
&=& f(0)+f'(0)x+ \int_0^x dt_1 \int_0^{t_1}dt_2~f''(0)+ \nonumber \\
&&~~~~~~~~~~~~~~~~~~~~~~~~~\int_0^x dt_1 \int_0^{t_1}dt_2~\int_0^{t_2}dt_3~f'''(t_3).
\eea

We can calculate that 
\be
\int_0^x dt_1 \int_0^{t_1}dt_2~f''(0) 
=  f''(0)\frac{1}{2}x^2.
\ee

Thus,
\be
f(x) = f(0)+f'(0)x+ f''(0)\frac{1}{2}x^2+ \int_0^x dt_1 \int_0^{t_1}dt_2\int_0^{t_2}dt_3~f'''(t_3).
\ee
This reaches the second-order Taylor polynomial. 

Further FTC applications can reach the 3rd-order and 4th-order, and higher-order Taylor polynomials. 
This derivation has a simple idea: Repeated use of FTC. The same OFFER principle has been applied to make the interior-boundary prediction
based on the simple FTC idea, and resulted in the theories for the conservation of energy, vorticity, and mass, represented by 
an integral along an open curve, over an open surface, and a solid body, corresponding to the gradient theorem, Stokes' theorem, and divergence theorem, respectively. 

\section{Rediscover the Displacement Current}
Displacement current is critically important for the discovery of electromagnetic waves. 
To rediscover the displacement current, we first ask a question. 

Inspired by Eq. (\ref{faraday1}), why don't we consider the following integral?
\be
\oint_C {\bf B}\cdot d{\bf r}
\label{eq1001}
\ee
where we treat ${\bf B}$ as a vector field. This integral is a magnetic circulation, in contrast to electric circulation in Eq. (\ref{faraday1}). 

\begin{figure} [h]
\centering
\includegraphics[width=4.0in]{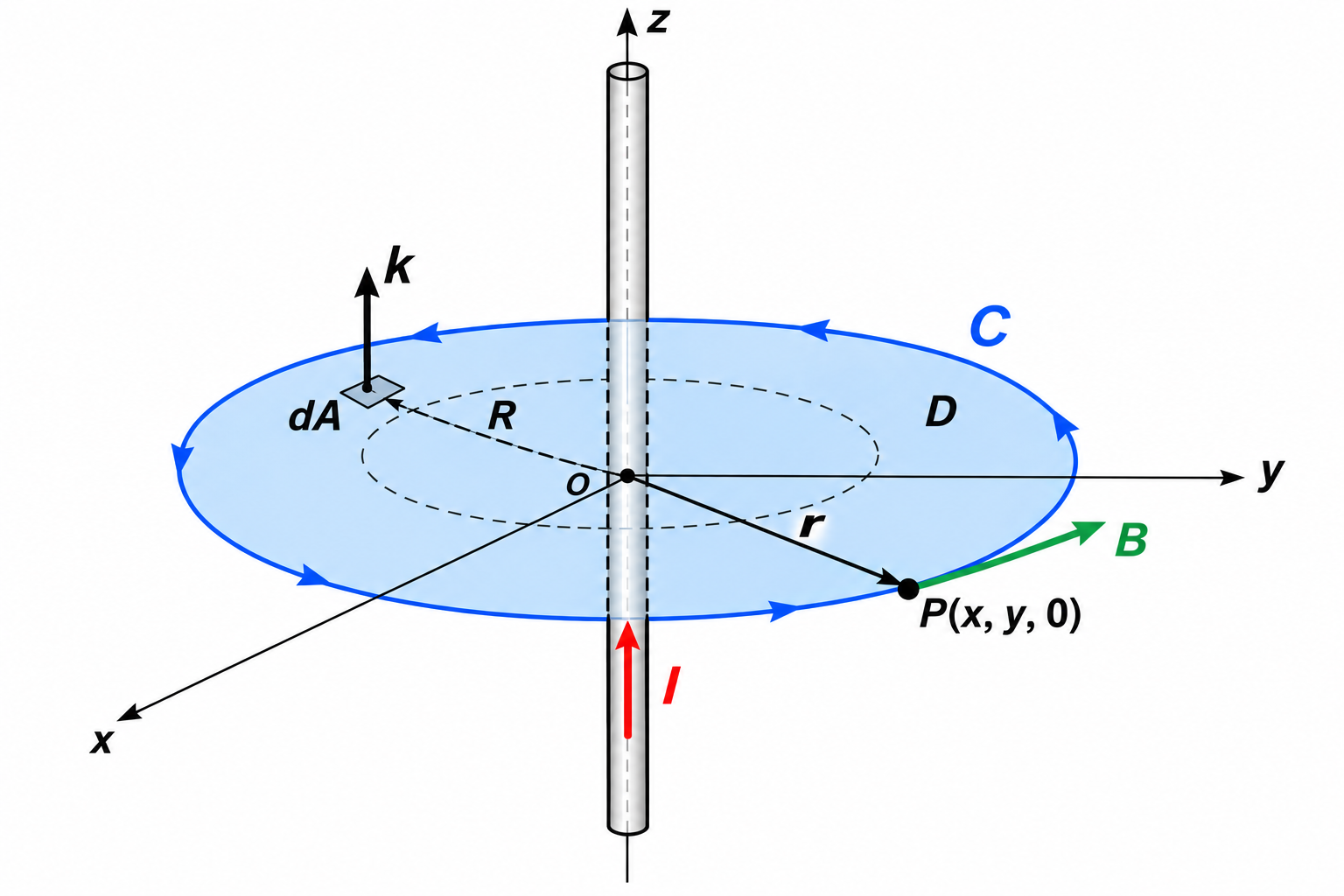}
\caption{Electric current $I$ and its induced magnetic field ${\bf B}$ in Amp\`ere's law. This figure was produced by ChatGPT AI
with a series of prompts proposed by the author. }
\label{fig10}
\end{figure}

Amp\`ere's law, named after French physicist and mathematician Andr\'e-Marie Amp\`ere (1775 – 1836), says that  an electric current $I$ 
can generate a rotational magnetic field ${\bf B}$ around the current. The quantitative relationship between $I$ and ${\bf B}$ 
can be described as follows 
\be
\oint_C {\bf B}\cdot d{\bf r} = \mu_0 I,
\label{eqamp1}
\ee
where the closed loop $C$ encloses the wire with the electric current $I$, and $\mu_0$ is a constant called the magnetic permeability. 
See Fig. \ref{fig10}.

Green's theorem implies 
\be
\oint_C {\bf B}\cdot d{\bf r} = \int_D (\nabla \times {\bf B})\cdot {\bf k} ~dA,
\label{eqamp2}
\ee
where $D$ is the region on the $xy-$plane enclosed by loop $C$, and ${\bf k}$ is the vertical unit vector along the $z$-direction. 

One can prove that 
\be
{\bf B} = \frac{\mu_0 I}{2\pi r^2} (-y, x),
\ee
and that 
\be
\nabla \times {\bf B} = \mu_0 I \delta (x, y)
\ee
where $r = |{\bf r}| = \sqrt{x^2 + y^2}$, and $\delta (x, y)$ is the Dirac delta function with 
\be
\int_D  \delta (x, y) ~dA = 1. 
\ee 
The technical details of these three formulas are beyond the first-time learners of calculus, but can be 
comprehended in the second or third semester of calculus. Nonetheless, the formulas provide insight for the 
calculus instructors. 

The above demonstrates the consistency between Green's theorem and Amp\`ere's law.

Now, we go back to our question at the beginning of this section:  What is the circulation defined by Eq. (\ref{eq1001})? 
A special case is that the wire with electric current $I$ inside the loop $C$ is absent, e.g., inside a vacuum. Green's theorem still implies
\be
\oint_C {\bf B}\cdot d{\bf r} = \int_D (\nabla \times {\bf B})\cdot {\bf k} ~dA.
\label{eqamp3}
\ee
What $\nabla \times {\bf B}$ should be when $I$ is zero? Inspired by Faraday's law
\be
\nabla \times {\bf E} = - {\bf B}_t
\label{eq1006}
\ee
and symbolic symmetry, we may guess that 
\be
\nabla \times {\bf B} = \mu_0 \epsilon_0 {\bf E}_t.
\label{eq1007}
\ee
Here, the constants $\mu_0 \epsilon_0$ help to balance the units (or dimensions) for both sides of the equation. 
Comparing this equation with Amp\`ere's law (\ref{eqamp1}), $\epsilon_0 {\bf E}_t$ must have the electric current dimension.  
This current is called the displacement current in Maxwell equations.   Scottish physicist and mathematician James Clerk Maxwell discovered this around 1861-1864. With or without the negative sign in Eqs. (\ref{eq1006}) and (\ref{eq1007}) are due to the right-hand rules. 

The equation below
\be
\nabla \times {\bf B} = \mu_0 \epsilon_0 {\bf E}_t + \mu_0 I
\label{eq1008}
\ee
is called the Amp\`ere-Maxwell law. Walking through the above process of Green's theorem, Amp\`ere's law, and inspired by Faraday's law, you have actually 
rediscovered the displacement current.

Another myth is the present form of Amp\`ere's law: electric current induces a circulating magnetic field. In fact, his discovery was about the 
interactive forces generated by two wires with electric current.  Danish physicist Hans Christian Ørsted's 1820's discovery is more directly 
linked to today's Amp\`ere's law: when current flowed through a wire, a nearby compass needle turned. 

 In terms of causality, we may say that the electric current in a conducting wire 
generates a magnetic circulation around the wire. This magnetic circulation can be observed easily with a home experiment kit
you can buy from Amazon or other stores. You can watch YouTube demos, such as the one by high school teacher Mr. Woodward {\tt https://www.youtube.com/watch?v=bEH7hCdhfhY}.

Ørsted's work inspired several French physicists and mathematicians for further explorations in the 1820s. 
Jean-Baptiste Biot and F\'elix Savart quantified the strength of the force acting on the compass needles. 
Pierre-Simon Laplace and Andr\'e-Marie Amp\`ere extensively explored the mathematics on the subject. 
Strictly speaking, the current form of Amp\`ere's law in textbooks was recasted by later writers.  

Recall that Green's theorem was also discovered in the 1820s when George Green studied electromagnetism. The displacement current 
was discovered in the early 1960s. 





The importance of the discovery of displacement current can be demonstrated in many ways. One is electromagnetic waves, and another is LC-oscillators. 

Taking derivative of Eq.  (\ref{eq1007}) with respect to time and substituting Eq. (\ref{eq1006}) into the result yield 
equations for electromagnetic waves with their propagation speed equal to 
\be
c = \sqrt{\mu_0 \epsilon_0},
\ee
which is the speed of light. 

The justification of the existence of displacement current in most textbooks is often made using a capacitor and Amp\`ere's law. Only 
the existence of a time-varying electric field can avoid the contradiction of a charged two-plate capacitor. Although 
$\epsilon_0 {\bf E}_t$ is called the displacement current, it is only a time-varying electric field, not a normal electric current. 
The displacement current is a critical concept in the LC-oscillation circuits in numerous industrial applications, e.g., radio tuning, signal generators, and radar systems.

\section{Summary and Discussion}

We have developed an OFFER (One-Formula-Facilitates-Every-Relation) theory for calculus. The theory is supported by two axioms: 
(i) The accumulation axiom that the total distance a car has been driven in a period of time is equal to the end odometer reading minus the beginning reading, and (ii) the level axiom that an uneven surface can always be leveled. 

Axiom (i) turns calculus into a method for interior-boundary prediction and is expressed as the FTC. 
Axiom (ii) provides mathematical rigor to the derivative calculation by a differential quotient established by Fermat, and is expressed as the MVT. The MVT also unifies Descartes' method of tangents with Fermat's. 

Extension of (i) in a time interval into the only three kinds of geometric objects in 3D space (an open curve, an open surface, and a solid body) leads to the conservation of energy, vorticity, and mass, which are 
respectively expressed by a line integral of the potential gradient on an open curve, a surface integral of the curl of a flow on an open surface, and a volume integral of the divergence of a flow in a solid body domain. The derivative integrations on their domains are predicted by their original functions' integrals on the domains' boundaries,
and vice versa. This statement is described by three fundamental theorems: the FTC for a curve, Green's theorem (or Stokes' theorem) for an open 
surface, and the divergence theorem for a body. 
These theorems lead to some fundamental principles of gravitation and electromagnetism: The integration of the gradient along an open curve 
is equal to the difference of the gravitational potentials; Faraday's law and Ampere's law can be described by Green's theorem or Stokes' theorem; and Gauss' laws of electricity and magnetism can be described by the divergence theorem. 

Therefore, FTC is the {\it one} formula facilitating every relationship (OFFER) in calculus, and has profound consequences in physics, engineering, and other application fields. 

With the OFFER theory, we hope to help both instructors and students focus on the fundamental ideas, methods, and purposes of calculus so that students are motivated, and instructors are excited about learning and teaching. We free students from memorizing tedious, meaningless formulas. Instead, students can make the most use of their common sense, critical thinking, and creativity. In the learning-by-doing process, students can draw figures by hand and with AI, create in-depth AI prompts on calculus, solve practical application problems with real-world data, and tell stories about the calculus theory, methods, and applications. Students can thus spend more time on idea creation and problem solving, such as generating AI prompts to formulate calculus problems, interpreting mathematical expressions, and summarizing solutions. 

At the end, we discuss how to use the materials in this article. The paper is written first for calculus instructors, second for students, and finally for calculus users, such as engineering researchers, with the hope that our materials can help modernize calculus education in the era of AI and the digital economy. Instructors may develop a big picture of calculus, understand the connections between calculus and physics (e.g., gravitational theory, thermodynamics, and electromagnetism), accelerate the pace of Cal I, II, and III series, and tell better stories, because they possess all the skills of differentiation, integration, and simple applications of FTC, Green's theorem, and the divergence theorem. They may use this article to argue for the modernization of calculus education and develop their lecture notes, computer code, AI prompts, and other teaching materials. 

This article is not meant to help first-year undergraduate students or high school students earn an A in a traditional calculus class.  The calculation details in some sections are beyond their comprehension. Nonetheless, this article can help student readers grasp the fundamental idea of calculus. They may qualitatively know that calculus is related to physics laws. Most students forget Green's theorem shortly after their final exam in traditional calculus. If students spend 5-6 hours on this article, they likely will not forget Green's theorem for the rest of their lives. 

For calculus users, our cut-and-add examples may help you generate novel calculus formulations for your research or job needs in practice. Our discovery examples may motivate you to create new calculus application fields and develop new methods for science, engineering, and other communities. 

\vskip 0.25cm

\noindent {\bf Acknowledgments}\\
This work was partially supported by the U.S. National Science Foundation under Grant IIS 2324008; the U.S. National Oceanographic and Atmospheric Administration under Grants NA22SEC4810016 and NA24NESX405C0006-T1-01.

\vskip 0.25cm

\noindent {\bf References}
\begin{description}
\item{[1]} Newton, I. (1736). The method of fluxions and infinite series; with its application to the geometry of curve-lines (Translated from the author’s Latin original; with “A perpetual comment upon the whole work” by John Colson). Henry Woodfall, London.
\item{[2]} Shen, S.S.P.,  Somerville, R.C.J.  (2019). Climate Mathematics: Theory and Applications, Cambridge University Press, Cambridge.
\item{[3]} Leibniz, G. (1684). Nova methodus pro maximis et minimis. {\it Acta Eruditorum} (in Latin). 3: 467–473.
\item{[4]} Struick, D. J. (1969). A Source Book in Mathematics, 1200-1800, Princeton University Press, Princeton.
\end{description}

\vskip 0.5cm

\noindent {\bf SAMUEL S.P. SHEN} is Distinguished Professor of Mathematics and Statistics at San Diego State University, and Visiting Research Mathematician at Scripps Institution of Oceanography, University of California – San Diego. Formerly, he was McCalla Professor of Mathematical and Statistical Sciences at the University of Alberta, Canada, and President of the Canadian Applied and Industrial Mathematics Society. He has held visiting positions at the NASA Goddard Space Flight Center, the NOAA Climate Prediction Center, and the University of Tokyo. Shen holds a B.Sc. degree in Engineering Mechanics and a Ph.D. degree in Applied Mathematics. 

\end{document}